\documentclass[preprint]{elsarticle}%review
\usepackage{lineno}
\journal{JOUARNAL's NAME}

\usepackage{slashbox}
\usepackage{amsmath}
\usepackage{amssymb}
\usepackage{amsthm}
\usepackage{graphicx}
\usepackage{epic,eepic,epsfig}
\usepackage{color}
\usepackage{subfigure}  % add 1.22
\usepackage{placeins}   % add 1.22
\usepackage{multirow}
\usepackage{epstopdf}
\usepackage{varwidth}
\usepackage{float}
\usepackage[table]{xcolor}
\usepackage{bm}
\usepackage[title]{appendix}
\usepackage{diagbox}
\usepackage{booktabs}
\usepackage[title]{appendix}

\numberwithin{equation}{section}

\newtheorem{theorem}{Theorem}[section]

\newtheorem{example}{Example}[section]
\newtheorem{lemma}{Lemma}[section]
\newtheorem{remark}{Remark}[section]

\newtheorem{assumption}{Assumption}[section]

\definecolor{tabclr}{cmyk}{0,0,1,0}

\allowdisplaybreaks[4]

\usepackage{caption2}
\usepackage{algorithm}
\usepackage{algorithmic}
\usepackage{datetime2}

\allowdisplaybreaks[4]

\makeatletter
\def\ps@pprintTitle{%
 \let\@oddhead\@empty
 \let\@evenhead\@empty
 \let\@oddfoot\@empty
 \let\@evenfoot\@empty
}
\makeatother

\begin{document}
\title{Nonsmooth data error estimates for exponential Runge--Kutta methods and applications to split exponential integrators}

\author[bjut]{Qiumei Huang}
\ead{qmhuang@bjut.edu.cn}

\author[ins]{Alexander Ostermann}
\ead{alexander.ostermann@uibk.ac.at}

%\author[bjut,ins]{Gangfan Zhong}
\author[bjut,ins]{Gangfan Zhong}
\ead{gfzhong@emails.bjut.edu.cn}

%\cortext[cor]{Corresponding author}
\address[bjut]{School of Mathematics, Statistics and Mechanics, Beijing University of Technology, 100124 Beijing, China}
 
\address[ins]{Department of Mathematics, University of Innsbruck, Technikerstr.~13, 6020 Innsbruck, Austria}

%\fntext[fn1]{\textbf{\color{blue} Version: 2025-06-02}}

%%%%% Begin Abstract %%%%%%%%%%%
\begin{abstract}
We derive error bounds for exponential Runge--Kutta discretizations of parabolic equations with nonsmooth initial data. Our analysis is carried out in a framework of abstract semilinear evolution equations with operators having non-dense domain. In particular, we investigate nonsmooth data error estimates for the Allen--Cahn and the Burgers’ equation. As an application, we apply these nonsmooth data error estimates to split exponential integrators and derive a convergence result in terms of the data. 
\end{abstract}

%%%%% end %%%%%%%%%%%

\maketitle
%
%\tableofcontents
%
%\linenumbers

%%%%%%%%%%%%%%%%%%%%%%%%%%%%%%
%%%%%%%%% ~~~~Introduction~~~~~~ %%%%%%%%%
 %%%%%%%%%%%%%%%%%%%%%%%%%%%%%%
\section{Introduction}
Exponential integrators are an important tool for solving linear and nonlinear evolution equations. By treating the linear term exactly and approximating the nonlinearity in an explicitly way, they are able to solve stiff problems in an accurate and efficient way. They have been employed in various large-scale computations, with applications to delay differential equations \cite{DaiHuang23:350}, phase field models \cite{HuangQiao24:116981}, and integro-differential equations \cite{OstermannSaedpanah23:1405}, among others.
The stiff convergence of exponential integrators has been widely studied. Convergence results are available for exponential Runge--Kutta methods \cite{HochbruckOstermann05:1069,HochbruckOstermann05:323,LuanOstermann13:3431}, exponential Rosenbrock methods \cite{HochbruckOstermann09:786} and exponential multistep methods \cite{HochbruckOstermann11:889}. For a comprehensive overview of exponential integrators, we refer the reader to \cite{HochbruckOstermann10:209}. Note, however, that the convergence analysis established in the aforementioned studies is always carried out in terms of the solution, i.e., sufficient temporal (and spatial) regularity of the solution is required. However, this demands certain compatibility conditions that are frequently impractical in real-world applications.

In this paper, we investigate nonsmooth data estimates for the semilinear parabolic problem
\begin{equation}\label{Eqn:intro-semilinear-probolems}
u^\prime  =Au +  f(t,u),\quad u(0)=u_0.
\end{equation}
We study the convergence of exponential Runge--Kutta discretizations of \eqref{Eqn:intro-semilinear-probolems} in terms of the prescribed data. We carry out this analysis in an abstract framework of analytic semigroups. For implicit and linearly implicit Runge--Kutta methods, a similar analysis has been carried out in \cite{LubichOstermann96:279,OstermannThalhammer00:167}. In these works, the operator $A:D(A)\subset X \to X$ was assumed to be densely defined (cf.~\cite{CrouzeixThomee87:359,ElliottLarsson92:603,Larsson92,LiMa21:23}). As a consequence, the choice $X=C(\overline{\Omega})$ was not possible there.

Here, we will use more general framework (see \cite{Lunardi95Book}) that allows us to consider also operators with non-dense domain. In this framework, we can also make a choice $X=C(\overline{\Omega})$ with the related maximum norm, which seems to be appropriate for many problems in applied mathematics. Moreover, nonlinearities (even nonlinearities involving derivatives) can easily be treated in this space. The use of the related interpolation spaces allows us to derive sharp error estimates and enables us to study a large class of parabolic problems.

As an application, we consider the convergence property of a split exponential Runge--Kutta method in which the operator $A$ is splitted into $A=A_1+A_2$. It is well known that operator splitting methods suffer from order reduction caused by incompatibilities between the nonlinearity and the boundary condition. Numerous efforts have been made to explain precisely this phenomenon \cite{FaouOstermann15:161} and to develop strategies to address it \cite{EinkemmerOstermann15:A1577,EinkemmerOstermann16:A3741}. Recently, a novel split exponential integrator was developed based on a second order explicit exponential Runge--Kutta method. The core of it lies in approximating the exponential-like $\varphi_k$ functions by means of $\varphi_k(\tau (A_1+A_2))\approx  k! \varphi_k(\tau A_1)\varphi_k (\tau A_2)$. We carry out a sound convergence analysis of such integrators in terms of the data.

The outline of the paper is as follows. In Section~\ref{Sec:framework}, we summarize the employed abstract framework and introduce the class of required interpolation spaces. 
In Section~\ref{Sec:linear}, error estimates for linear parabolic problems in terms of data are studied. %based on assumptions regarding the bounds of the source function and its derivatives.
In Section~\ref{Sec:semilinear},  we analyze nonsmooth data error estimates for semilinear parabolic problems. We begin by deriving an estimate for the derivative of the solution, which leads to a first order convergence result under the assumption of local Lipschitz continuity for the nonlinearity. If the initial data satisfies a compatibility condition and the nonlinearity is smoother, a second order convergence result can be achieved. In Section~\ref{Sec:example}, we study nonsmooth data error estimates for both the Allen–Cahn and the Burgers’ equation. In Section~\ref{Sec:split}, we apply the non\-smooth data error estimates to a split exponential integrator and derive a convergence result in terms of the data.

\section{Analytical framework}\label{Sec:framework}
Our analysis below will be carried out in a framework of abstract semilinear parabolic evolution equations with non-dense domain, as developed in \cite[Chapter 7]{Lunardi95Book}. Let $X$ be a Banach space, $A$ an operator on $X$ and let $D(A)$ denote the domain of $A$ in $X$.  In contrast to standard assumptions as presented in \cite{EngelNagel00Book,Henry81Book,Pazy12Book}, the operator $A$ is not assumed to be densely defined in $X$, i.e., we have only $\overline{D(A)}\subset X$. This enables us to consider $X=C(\overline{\Omega})$, where the associated norm is particularly suitable for practical use, as it provides a direct measure of the pointwise numerical error commonly encountered in simulations. Our main assumption is the following one.

\begin{assumption}\label{Ass:sectorial-A}
Let $X$ be a Banach space with norm $\|\cdot\|$. We assume that the linear operator $A:D(A)\subset X \to X$ is the infinitesimal generator of an analytic semigroup $\mathrm{e}^{tA}$ on $X$.
\end{assumption}

Under this assumption, a class of real interpolation spaces between $X$ and $D(A)$ is given as (see \cite{Lunardi95Book})
$$
D_A(\alpha, p)= \big\{x \in X: t \mapsto v(t)= \|t^{1-\alpha-1 / p} A \mathrm{e}^{t A} x \| \in L^p(0,1) \big\},\quad 0<\alpha <1,~1\leq p \leq \infty
$$
with corresponding norm
$$
\|x\|_{D_A(\alpha, p)}=\|x\|+\|v\|_{L^p(0,1)}.
$$
For convenience, we set $D_A(0,p)=X$ for $p\in [1,\infty]$. For $0<\alpha <1$, let $X_\alpha$ be any Banach space satisfying
$$
D_A(\alpha,1)\subset X_\alpha \subset D_A(\alpha,\infty)
$$
with continuous embeddings.
The norm of the space $X_\alpha$ is denoted as $\|\cdot\|_\alpha$. Moreover, we set  $X_0=X$. We further make the following assumption.

\begin{assumption}\label{Ass:sectorial-Aalpha}
The part of $A$ in $X_\alpha$ is sectorial in $X_\alpha$.
\end{assumption}

Note that this assumption is satisfied if $X_\alpha$ is the real interpolation space $D_A(\alpha,p)$ or the domain $D((-A + \omega I)^\alpha)$ of a fractional power of $-A + \omega I$ (see \cite[p.~253]{Lunardi95Book}), where $\omega$ is a real number such that $\text{Re}\, \sigma(-A+\omega I) >0$, where $\sigma$ denotes spectrum. We recall that under the above assumptions the operator $A$ satisfies the properties
\begin{equation}\label{Est:smoothing-etA}
\|\mathrm{e}^{t A}\|+\|tA\mathrm{e}^{t A}\| + \|\mathrm{e}^{t A}\|_{\alpha}   \leq C, \quad  0\leq t \leq T.
\end{equation}

\section{Linear problem}\label{Sec:linear} 
In this section, we consider linear parabolic problems 
\begin{equation}\label{Eqn:linear-probolems}
u^\prime  =Au + g(t),\quad u(0)=u_0.
\end{equation}
Applying an $s$-stage exponential Runge--Kutta method \cite{HochbruckOstermann05:1069} to this linear problem with a time step $\tau>0$ yields at $t_n=n \tau$ an approximation $u_n$ to $u(t_n)$, which is recursively given by 
\begin{align}
u_{n+1} & = \mathrm{e}^{ \tau A}u_n+\tau \sum_{i=1}^s b_i(\tau A)g(t_{ni}),  
\quad t_{ni}=t_n+c_i\tau, \label{ERK-linear} 
\end{align}
where $c_i \in [0,1]$  for $1\leq i \leq s$. The weights $b_i(z)$ are linear combinations of the functions $\varphi_k(z)$, defined by
\begin{equation}\label{Eqn:varphik}
\varphi_0(z)=\mathrm{e}^z,\quad \varphi_k(z)=\int_0^1 \mathrm{e}^{(1-\xi)z}\frac{\xi^{k-1}}{(k-1)!}\,\mathrm{d}\xi,\quad k\geq1.
\end{equation}
The following estimate is frequently used throughout the analysis,
\begin{equation}\label{est:X-DAalpha1}
\|t^{\beta} \mathrm{e}^{t A}\|_{L(X, D_A(\beta, 1))}\leq C, \quad 0\leq\beta<1,~0<t\leq T.
\end{equation}
We recall its proof from \cite[Proposition 2.2.9]{Lunardi95Book} below to make the reader familiar with the used arguments in the proof. Using the defining property of interpolation spaces and the bounds \eqref{Est:smoothing-etA}, we get the bound
$$
 \|t A \mathrm{e}^{t A} x \|_{D_A(\beta, 1)} \leq C \|t  A  \mathrm{e}^{t A} x \|_{D(A)}^\beta \|t  A  \mathrm{e}^{t A} x \|^{1-\beta} \leq Ct^{-\beta}\|x\|, \quad 0 < t \leq 1.
$$
Further, we have for $0 < t \leq 1$,
$$
\begin{aligned}
  \|\mathrm{e}^{t A} x \|_{D_A(\beta, 1)}& \leq \|\mathrm{e}^A x \|_{D_A(\beta, 1)}+ \bigg\|\int_t^1 A \mathrm{e}^{\xi A} x \,\mathrm{d} \xi  \bigg\|_{D_A(\beta, 1)} \\
  &  \leq  C  \|\mathrm{e}^{A}x\|_{D(A)}^\beta \|\mathrm{e}^{A}x\|^{1-\beta}+C  \int_t^1 \xi^{-\beta-1} \|x\| \,\mathrm{d}\xi  \leq Ct^{-\beta}\|x\| ,
\end{aligned}
$$
which leads to \eqref{est:X-DAalpha1} for $0<t\leq 1$.
For $1<t \leq T$, we have
$$
\begin{aligned}
\| \mathrm{e}^{t A} \|_{L (X, D_A(\beta, 1) )}   \leq \|\mathrm{e}^{  A} \|_{L (X, D_A(\beta, 1) )}\|\mathrm{e}^{(t-1) A} \| \leq C .
\end{aligned}
$$
The combination of the above inequalities leads to \eqref{est:X-DAalpha1}. Due to \eqref{Eqn:varphik}, we also have 
\begin{equation}\label{Est:varphik-X-DAalpha1}
\| t^\beta   \varphi_k( t A) \|_{L(X, D_A(\beta, 1))}  \leq C,\quad 0\leq  \beta < 1,~0<t \leq T.
\end{equation}

Let $g:[0,T]\to X$. If $g\in L^\infty( (0,T) ;X)$, it follows from \cite[Proposition 4.2.1]{Lunardi95Book} that 
$$
\int_0^t \mathrm{e}^{(t-\xi)A} g(\xi) \, \mathrm{d}\xi \in C^{1-\gamma} ( [0,T] ;D_A(\gamma,1)) \quad \mbox{for any $\gamma \in (0,1)$}.
$$
Given $u_0\in X$, the function $u$ defined by
\begin{equation}\label{Eqn:g-linear-mild}
u(t)=\mathrm{e}^{t A}u_0 + \int_0^t \mathrm{e}^{(t-\xi)A} g(\xi)\,\mathrm{d}\xi,\quad 0\leq t \leq T,
\end{equation}
is called the mild solution of \eqref{Eqn:linear-probolems}.

We begin with error estimates for first order methods. Such methods satisfy the order condition (see \cite{HochbruckOstermann05:1069,LuanOstermann13:3431})
\begin{equation}\label{Eqn:linear-first-order}
 \sum_{i=1}^s b_i(z) = \varphi_1(z).
\end{equation}

\begin{theorem}\label{Thm:linear-first}
Let the linear initial value problem \eqref{Eqn:linear-probolems} satisfy Assumption \ref{Ass:sectorial-A} and
\begin{equation}\label{Eqn:g-linear}
\|g(t)\|\leq B_1, \quad 0\leq t \leq T,\qquad \|g^\prime(t)\| \leq B_1t^{-1},\quad 0< t \leq T.
\end{equation}
Consider for its numerical solution an exponential Runge--Kutta method \eqref{ERK-linear} satisfying \eqref{Eqn:linear-first-order}. For initial data $u_0\in X_\alpha$, the error satisfies the bound
$$
\|u_n -u(t_n)\|_{\alpha} \leq  C\tau^{1-\alpha}( \left|\log \tau\right| +1)
$$
uniformly in $0\leq t_n \leq T$. The constant $C$ depends on $B_1$ and $T$, but is independent of $n$ and~$\tau$.
\end{theorem}
\begin{proof}
Solving the recursion \eqref{ERK-linear} gives
\begin{equation}\label{Eqn:linear-un-recursion}
u_n = \mathrm{e}^{n\tau A}u_0 + \tau \sum_{k=0}^{n-1}\mathrm{e}^{(n-k-1)\tau A }\sum_{i=1}^s b_i(\tau A)g(t_{ki}).
\end{equation}
We rewrite the mild solution \eqref{Eqn:g-linear-mild} at $t=t_n$ as
\begin{equation}\label{Eqn:linear-utn-recursion}
u(t_n)=\mathrm{e}^{n \tau A}u_0 + \sum_{k=0}^{n-1} \mathrm{e}^{(n-k-1)\tau A} \int_0^{\tau} \mathrm{e}^{(\tau-\xi)A} g(t_k+\xi)\,\mathrm{d}\xi.
\end{equation}
Subtracting \eqref{Eqn:linear-utn-recursion} from \eqref{Eqn:linear-un-recursion} yields
\begin{equation}\label{Eqn:error-exposition}
u_n -u(t_n) =   \tau \sum_{k=0}^{n-1} \mathrm{e}^{(n-k-1)\tau A} \bigg( \sum_{i=1}^s b_i(\tau A)g(t_{ki}) - \frac{1}{\tau}\int_0^{\tau} \mathrm{e}^{(\tau-\xi)A} g(t_k+\xi)\,\mathrm{d}\xi \bigg).
\end{equation}
For $n=1$, using the boundedness of $g$ in $[0,T]$ and the estimates \eqref{est:X-DAalpha1} and \eqref{Est:varphik-X-DAalpha1}, we have
\begin{align}
\|u_1 -u(\tau) \|_{\alpha} &\leq C\|u_1 -u(\tau) \|_{D_A(\alpha,1)} \notag \\
&\leq C\tau^{1-\alpha}\sum_{i=1}^s \|\tau^{\alpha} b_i(\tau A)\|_{L(X,D_A(\alpha,1))}  \|g(t_{0i})\|  \notag \\
&\quad+ C\int_0^{\tau} (\tau-\xi)^{-\alpha}\|(\tau-\xi)^{\alpha} \mathrm{e}^{(\tau-\xi)A}\|_{L(X,D_A(\alpha,1))}  \|g(\xi) \|  \,\mathrm{d}\xi \notag \\
& \leq C\tau^{1-\alpha}.\label{Eqn:linear-n=1}
\end{align}
We now consider the case of $n\geq 2$. With the help of the Taylor expansions 
$$
\begin{aligned}
g(t_{ni})   = g(t_n) + \int_0^{c_i\tau} g^\prime ( t_n+\xi)\,\mathrm{d}\xi, \quad
g(t_{n}+\xi )   = g(t_n) + \int_0^{\xi} g^\prime ( t_n+\eta)\,\mathrm{d}\eta 
\end{aligned}
$$
and the use of \eqref{Eqn:linear-first-order}, we have
\begin{align}
\tau\sum_{i=1}^s b_i(\tau A)g(t_{ni}) & = \tau \varphi_1(\tau A) g(t_{n})+\delta_{n1}~ \mbox{with}~\delta_{n1}=\tau \sum_{i=1}^s b_i(\tau A)\int_0^{c_i\tau} g^\prime ( t_n+\xi)\,\mathrm{d}\xi, \label{Eqn:defectn1} \\
\int_0^{\tau} \mathrm{e}^{(\tau-\xi)A} g(t_n+\xi)\,\mathrm{d}\xi  & = \tau \varphi_1(\tau A) g(t_{n})+\widetilde{\delta}_{n1}~ \mbox{with}~\widetilde{\delta}_{n1}=\int_0^{\tau} \mathrm{e}^{(\tau-\xi)A}\int_0^{\xi} g^\prime ( t_n+\eta)\,\mathrm{d}\eta \, \mathrm{d}\xi.\label{Eqn:tilde-defectn1}
\end{align}
Using \eqref{Eqn:g-linear} and the estimates \eqref{est:X-DAalpha1} and \eqref{Est:varphik-X-DAalpha1}, the remainders are bounded by
\begin{align}
\|\delta_{n1}\|_{\alpha},\, \|\widetilde{\delta}_{n1}\|_{\alpha} & \leq   C t_n^{-1}\tau^{2-\alpha} .\label{Eqn:defectn1-bound}
\end{align}
Plugging \eqref{Eqn:defectn1}-\eqref{Eqn:tilde-defectn1} into \eqref{Eqn:error-exposition} yields
\begin{equation}\label{Eqn:error-exposition-rewrite}
u_n - u(t_n) =   \sum_{k=0}^{n-1} \mathrm{e}^{(n-k-1)\tau A} ( \delta_{k1}-\widetilde{\delta}_{k1}). 
\end{equation}
Using \eqref{Eqn:defectn1-bound}, we finally estimate \eqref{Eqn:error-exposition-rewrite} by
$$
\|u_n -u(t_n)\|_{\alpha} \leq C\tau \sum_{k=1}^{n-1} \tau^{1-\alpha}t_k^{-1}+C\tau^{1-\alpha} \leq C\tau^{1-\alpha}( \left|\log \tau\right| +1),
$$
which is the desired result.
\end{proof}

Next, we consider the error estimates for second order methods. These methods satisfy the order conditions (see \cite{HochbruckOstermann05:1069,LuanOstermann13:3431})
\begin{equation}\label{Eqn:linear-second-order}
 \sum_{i=1}^s b_i(z) = \varphi_1(z),\quad \sum_{i=1}^s b_i(z)c_i = \varphi_2(z).
 \end{equation}

\begin{theorem}\label{Thm:linear-second}
Let the linear initial value problem \eqref{Eqn:linear-probolems} satisfy Assumption \ref{Ass:sectorial-A} and
\begin{equation}\label{Eqn:g-linear-second}
\|g(t)\|,\, \|g^\prime(t)\|\leq B_2, \quad 0\leq t \leq T,\qquad \|g^{\prime\prime}(t)\| \leq B_2 t^{-1},\quad 0< t \leq T.
\end{equation}
Consider for its numerical solution an exponential Runge--Kutta method \eqref{ERK-linear} satisfying \eqref{Eqn:linear-second-order}. For initial data $u_0\in X_\alpha$, the error satisfies the bound
$$
\|u_n -u(t_n)\|_{\alpha} \leq  C\tau^{2-\alpha}( \left|\log \tau \right| +1)
$$
uniformly in $0\leq t_n  \leq T$. The constant $C$ depends on $B_2$ and $T$, but is independent of $n$ and~$\tau$.
\end{theorem}
\begin{proof}
The proof is similar to that of Theorem \ref{Thm:linear-first}. For $n=1$, using \eqref{Eqn:defectn1}-\eqref{Eqn:tilde-defectn1} and \eqref{Eqn:g-linear-second}, we derive
$$
\|u_1-u(\tau)\|_{\alpha} =     \Bigg\| \tau \sum_{i=1}^s b_i(\tau A)\int_0^{c_i\tau} g^\prime ( \xi)\,\mathrm{d}\xi -\int_0^{\tau} \mathrm{e}^{(\tau-\xi)A}\int_0^{\xi} g^\prime ( \eta)\,\mathrm{d}\eta\,\mathrm{d}\xi      \Bigg\|_{\alpha}
\leq C\tau^{2-\alpha}.
$$
For $n \geq 2$, with the help of the Taylor expansions
$$
\begin{aligned}
g(t_{ni}) & = g(t_n) + c_i g^\prime(t_n)+ \int_0^{c_i\tau} (c_i \tau -\xi)g^{\prime\prime} ( t_n+\xi)\,\mathrm{d}\xi, \\
g(t_{n}+\xi ) & = g(t_n) +\xi g^\prime(t_n)+ \int_0^{\xi} (\xi - \eta) g^{\prime\prime} ( t_n+\eta)\,\mathrm{d}\eta, \\
\end{aligned}
$$
and the use of \eqref{Eqn:linear-second-order}, we have
\begin{align}
\tau\sum_{i=1}^s b_i(\tau A)g(t_{ni}) & = \tau \varphi_1(\tau A) g(t_{n})+\tau^2 \varphi_1(\tau A) g^\prime(t_{n})+\delta_{n2} \quad \mbox{with}\quad\|\delta_{n2}\|_{\alpha} \leq Ct_n^{-1} \tau^{3-\alpha}, \notag  \\
\int_0^{\tau} \mathrm{e}^{(\tau-\xi)A} g(t_j+\xi)\,\mathrm{d}\xi  & = \tau \varphi_1(\tau A) g(t_{n})+\tau^2 \varphi_1(\tau A) g^\prime(t_{n})+\widetilde{\delta}_{n2}
\quad \mbox{with}\quad\|\widetilde{\delta}_{n2}\|_{\alpha} \leq Ct_n^{-1} \tau^{3-\alpha}. \notag
\end{align}
Plugging these inequalities into \eqref{Eqn:error-exposition} and estimating the resulting term as before leads to the desired result.
\end{proof}

\section{Semilinear problem}\label{Sec:semilinear}

In this section, we will derive error bounds for explicit exponential Runge–Kutta discretizations of semilinear parabolic problems with nonsmooth initial data
\begin{equation}\label{Eqn:semilinear-probolems}
u^\prime  =Au +  f(t,u),\quad u(0)=u_0.
\end{equation}
Applying an $s$-stage explicit exponential Runge--Kutta (EERK) method to this problem yields
\begin{subequations}\label{ERK-semilinear}
\begin{align}
    u_{n+1} &= \mathrm{e}^{\tau A} u_n + \tau \sum_{i=1}^{s} b_i (\tau A) F_{ni}, \label{ERK-semilinear-1} \\
    U_{ni} &= \mathrm{e}^{c_i\tau A} u_n + \tau \sum_{j=1}^{i-1} a_{ij} (\tau A) F_{nj},\quad i=1,\ldots,s, \label{ERK-semilinear-2} \\
    F_{nj} &= f(t_{nj}, U_{nj}),\quad j=1,\ldots,s,  \label{ERK-semilinear-3} 
\end{align}
\end{subequations}
where $c_i \in [0,1]$  for $1\leq i \leq s$ and $c_1=0$. The weights $a_{ij}(z)$ and $b_i(z)$ are linear combinations of the functions $\varphi_k(z)$. We recall some estimates for $\mathrm{e}^{tA}$, which can be found in the proof of \cite[Proposition 2.2.9]{Lunardi95Book}.

\begin{lemma}\label{Lem:est-etA}
Under Assumption \ref{Ass:sectorial-A}, the following estimates hold uniformly for $0<t\leq T$:
\begin{itemize}\rm
\item[(a)] $\|t^{-\beta_1+\beta_2} \mathrm{e}^{tA} \| _{L(D_A(\beta_1,\infty),D_A(\beta_2,\infty))} \leq C,\quad 0\leq \beta_1\leq \beta_2<1$,
\item[(b)] $\|t^{1-\beta_1+\beta_2}A\mathrm{e}^{tA}\|_{L(D_A(\beta_1,\infty),D_{A}(\beta_2,1))} \leq C,\quad 0\leq \beta_1,\beta_2 <1$.
%\item[(d)] $\|t^{1-\beta_1}A\mathrm{e}^{tA}\|_{L(D_A(\beta_1,\infty),X)}\leq C,\quad 0\leq \beta_1 <1$.
\end{itemize}
\end{lemma}

Our basic assumption on the nonlinearity $f$ is the following. 
\begin{assumption}\label{Ass:f-lip}
Let $f:[0,T]\times X_\alpha \to X$ be locally Lipschitz bounded for some $\alpha \in  [0,1)$: for every $R>0$, there exists $L=L(R,T)$ such that
\begin{equation*}
\|f(t_1,v_1)-f(t_2,v_2)\| \leq L \big( |t_1-t_2| + \|v_1-v_2\|_{\alpha} \big)
\end{equation*}
for all $t\in[0,T]$ and $v_1,v_2\in X_\alpha$ with $\max\big( \|v_1\|_{\alpha},\|v_2\|_{\alpha})\leq R$.
\end{assumption}

%A function $u \in C((0,T];X_\alpha)$ is called a mild solution of problem \eqref{Eqn:semilinear-probolems} in the interval $[0,T]$ if it is a solution of the integral equation
%\begin{equation}\label{Eqn:mild-semilinear}
%u(t)=\mathrm{e}^{t A}u_0 + \int_0^t \mathrm{e}^{(t-\xi)A} f(\xi,u(\xi))\,\mathrm{d}\xi,\quad 0\leq t \leq T 
%\end{equation}
%and such that $f(\cdot,u(\cdot))$ belongs to $L^1( (0,T) ;X)$.
Based on Assumption \ref{Ass:f-lip}, we obtain from  \cite[Theorem 7.1.2]{Lunardi95Book} the following result: Given any initial data $u_0\in X_\alpha$, problem \eqref{Eqn:semilinear-probolems} has a unique mild solution $u \in L^\infty( (0,T);X_\alpha)\cap C((0,T];X_\alpha)$ for some $T>0$. For the sake of convenience, we define
$$
g(t)=f(t,u(t)).
$$ 
With the help of \eqref{Est:smoothing-etA} and \eqref{est:X-DAalpha1}, the solution is shown to be bounded as follows:
\begin{equation*}
\|u(t)\|_{\alpha} \leq C\|u_0\|_{\alpha} + CT^{1-\alpha}  \|g\|_{L^\infty(0,T;X)}, \quad 0\leq t \leq T .
\end{equation*}

To estimate the derivative of the solution, we extend \cite[Theorem 3.5.2]{Henry81Book} to the case of problems with non-dense domain, which necessitates appropriate modifications of the original proof.

\begin{lemma}\label{Lem:derivative}
Let the initial value problem \eqref{Eqn:semilinear-probolems} with $u_0\in X_\alpha\cap D_A(\gamma,\infty)$ for some $\gamma \in [\alpha ,1)$ satisfy Assumptions \ref{Ass:sectorial-A}-\ref{Ass:sectorial-Aalpha} and \ref{Ass:f-lip}. Then the derivative of the solution in $D_A(\mu,1)$ for $\mu\in [0,1)$ is bounded by
$$
\|u^\prime(t)\|_{D_A(\mu,1)} \leq  Ct^{\gamma-\mu-1},\quad 0< t \leq T,
$$
where the constant $C$ depends on $T$. 
\end{lemma}
\begin{proof}
Let $\alpha \leq \gamma< \beta<1$, $0<\sigma<T/2$ and recall $g(t)=f(t,u(t))$. It holds 
$$
\begin{aligned}
\| u(\sigma)\|_{D_A(\beta,\infty)} & \leq  \| \mathrm{e}^{\sigma A}u_0\|_{D_A(\beta,\infty)} + C\int_0^\sigma \|\mathrm{e}^{(\sigma-\xi)A}g(\xi)\|_{D_A(\beta,1)}\,\mathrm{d}\xi
\\
& \leq C\sigma^{\gamma-\beta}\|u_0\|_{D_A(\gamma,\infty)}+C\int_0^\sigma (\sigma-\xi)^{-\beta}\|g(\xi)\|  \,\mathrm{d}\xi\leq C\sigma^{\gamma-\beta},
\end{aligned}
$$
where property (a) of Lemma \ref{Lem:est-etA} and \eqref{est:X-DAalpha1} were used.
%$$
%\| u(\sigma)\|_{D_A(\beta,\infty)} \leq C\sigma^{\alpha-\beta}\|u_0\|_{D_A(\alpha,\infty)}+\int_0^\sigma (\sigma-s)^{-\beta}\|g(s)\|  \,\mathrm{d}s\leq C\sigma^{\alpha-\beta}.
%$$
For $\sigma<t<t+h\leq T$, noting that
$$
\begin{aligned}
u(t+h) - u(t)& = \mathrm{e}^{(t+h-\sigma)A}u(\sigma) +\int_\sigma^{\sigma+h} \mathrm{e}^{(t+h-\xi)A}g(\xi)\,\mathrm{d}\xi +\int_{\sigma+h}^{t+h} \mathrm{e}^{(t+h-\xi)A}g(\xi)\,\mathrm{d}\xi \\
&\quad -\mathrm{e}^{(t-\sigma)A}u(\sigma)-\int_\sigma^{t} \mathrm{e}^{(t-\xi)A}g(\xi)\,\mathrm{d}\xi 
\\
&=(\mathrm{e}^{hA}-1)\mathrm{e}^{(t-\sigma)A}u(\sigma)+\int_{\sigma}^{\sigma+h} \mathrm{e}^{(t+h-\xi)A}g(\xi)\,\mathrm{d}\xi\\
&\quad +\int_{\sigma}^t \mathrm{e}^{(t-\xi)A}\big(g(\xi+h) - g(\xi) \big)\, \mathrm{d}\xi ,
\end{aligned}
$$
we obtain from Assumption \ref{Ass:f-lip}
\begin{align}
&\|g(t+h) - g(t)\| \notag \\
& \quad\leq Lh + L\|u(t+h)-u(t)\|_{\alpha} \notag \\
& \quad\leq Lh + CL\int_0^{h} \|A\mathrm{e}^{\xi A} \mathrm{e}^{(t-\sigma)A}u(\sigma)\|_{D_A(\alpha,1)} \,\mathrm{d}\xi +CL\int_{\sigma}^{\sigma+h} \|\mathrm{e}^{(t+h-\xi)A}g(\xi)\|_{D_A(\alpha,1)} \,\mathrm{d}\xi \notag
\\
&\quad\quad +CL\int_{\sigma}^t \big\| \mathrm{e}^{(t-\xi)A}\big(g(\xi+h) - g(\xi) \big)\big\|_{D_A(\alpha,1)} \,\mathrm{d}\xi  \notag \\
&\quad \leq  Lh +CL \int_0^{h} (t-\sigma+\xi)^{\beta-\alpha-1}\|u(\sigma)\|_{D_A(\beta,\infty)} \,\mathrm{d}\xi +CL\int_{\sigma}^{\sigma+h} (t+h-\xi)^{-\alpha}\|g(\xi)\|  \,\mathrm{d}\xi
\notag \\
&\quad\quad +CL\int_{\sigma}^t (t-\xi)^{-\alpha}\|  g(\xi+h) - g(\xi)  \|  \,\mathrm{d}\xi   \notag \\
&\quad\leq  Lh + CL h\big( (t-\sigma)^{\beta-\alpha-1}\|u(\sigma)\|_{D_A(\beta,\infty)}+(t-\sigma)^{-\alpha}\big)   +CL\int_{\sigma}^t (t-\xi)^{-\alpha}  \|g(\xi+h) - g(\xi)  \|  \,\mathrm{d}\xi,
\notag
\end{align}
where property (b) of Lemma \ref{Lem:est-etA} and \eqref{est:X-DAalpha1} were used.
An application of Gronwall's lemma \cite[Lemma 6.3]{ElliottLarsson92:603} yields
\begin{equation}\label{Eqn:g-lipschitz}
\|g(t+h) - g(t)\| \leq Ch\big( (t-\sigma)^{\beta-\alpha-1}\|u(\sigma)\|_{D_A(\beta,\infty)} +(t-\sigma)^{-\alpha}\big).
\end{equation}

Note that
$$
\begin{aligned}
u(t) & = \mathrm{e}^{tA} u_0 + \int_0^t \mathrm{e}^{(t-\xi)A}g(\xi)\, \mathrm{d} \xi \\
& =  \mathrm{e}^{tA} u_0 +  \int_0^t \mathrm{e}^{(t-\xi)A}\big( g(\xi) - g(t)\big) \, \mathrm{d} \xi  + \int_0^t \mathrm{e}^{(t-\xi)A}g(t)\, \mathrm{d} \xi
\end{aligned}
$$
and $A\mathrm{e}^{(t-\cdot)A} \big( g(\cdot) - g(t)\big)\in L^1( (0,t); X)$ due to \eqref{Eqn:g-lipschitz}.
The derivative of the solution for $t>0$ is given as
$$
\begin{aligned}
u^\prime(t) & = A\mathrm{e}^{tA} u_0 +  A\int_0^t \mathrm{e}^{(t-\xi)A}\big( g(\xi) - g(t)\big) \, \mathrm{d} \xi + A\int_0^t \mathrm{e}^{(t-\xi)A}g(t)\, \mathrm{d} \xi  + g(t)  \\
 & =  A\mathrm{e}^{tA} u_0 +  \int_0^t A\mathrm{e}^{(t-\xi)A} \big( g(\xi) - g(t)\big) \, \mathrm{d} \xi + \mathrm{e}^{tA} g(t) ,
\end{aligned}
$$
for which \cite[(ii) of Proposition 2.1.4]{Lunardi95Book} was used.

Finally, by using property (b) of Lemma \ref{Lem:est-etA} and \eqref{est:X-DAalpha1}, the derivative of the solution is bounded by
\begin{align}
\|u^\prime(t)\|_{_{D_A(\mu,1)}} 
&\leq  Ct^{\gamma-\mu-1} \|u_0\|_{D_A(\gamma,\infty)}+ C\int_0^{2\sigma} (t-\xi)^{-\mu-1}\| g(\xi) -g(t) \|\,\mathrm{d}\xi \notag \\
&\quad+C \int_{2\sigma}^t (t-\xi)^{-\mu-1}\| g(\xi) -g(t) \|\,\mathrm{d}\xi + Ct^{-\mu}\|g\|_{L^\infty(0,T;X)}.\notag
\end{align}
With the help of \eqref{Eqn:g-lipschitz}, we obtain
\begin{align}
\|u^\prime(t)\|_{_{D_A(\mu,1)}} 
&\leq Ct^{\gamma-\mu-1} \|u_0\|_{D_A(\gamma,\infty)} +C\big( (t-2\sigma)^{-\mu} -t^{-\mu}\big) \notag \\
&\quad+ C\int_{2\sigma}^t (t-\xi)^{-\mu} \big( (\xi-\sigma)^{\beta-\alpha-1}\|u(\sigma)\|_{D_A(\beta,\infty)} +(\xi-\sigma)^{-\alpha}\big) \,\mathrm{d}\xi + Ct^{-\mu}\|g\|_{L^\infty(0,T;X)}   
\notag \\
&\leq C\big(t^{\gamma-\mu-1}+ (t-2\sigma)^{-\mu}+ t^{-\mu} +(t-2\sigma)^{1-\mu}(\sigma^{\gamma-\alpha-1} + \sigma^{-\alpha})    \big).\notag
\end{align}
The proof is completed by setting $\sigma=t/4$. 
\end{proof}

To carry out the convergence analysis, we follow the ideas in \cite{LubichOstermann96:279} and decompose $u$ into two parts as $u=x+y$, where $x$ and $y$ are defined by the following two parabolic systems:
\begin{subequations}\label{u=x+y} 
\begin{align}
x^\prime& = Ax,\phantom{+ g(t) }  ~ \quad x(0)=u_0,  \label{Eqn:x}\\
y^\prime& = Ay + g(t), \quad y(0)=0, \label{Eqn:y}
\end{align}
\end{subequations}
where $g(t)=f(t,u(t))$. Note that \eqref{Eqn:y} is a linear parabolic problem. Therefore, the results from Section~\ref{Sec:linear} are appliable. This is the reason why we perform the above decomposition.

The solution of \eqref{Eqn:x} is simply given by
\begin{align}
x(t)&=\mathrm{e}^{t A} u_0. \label{solution-x}
\end{align}
Applying an EERK method to \eqref{Eqn:y} gives the approximations $y_{n+1}\approx y(t_{n+1})$ and $Y_{ni}\approx y(t_{ni})$ as follows:
\begin{align}
    y_{n+1} &= \mathrm{e}^{\tau A} y_n + \tau \sum_{i=1}^{s} b_i (\tau A) g(t_{ni}), \label{ERK-y-1}
      \\
    Y_{ni} &= \mathrm{e}^{c_i\tau A} y_n + \tau \sum_{j=1}^{i-1} a_{ij} (\tau A)g(t_{nj}),\quad i=1,\ldots,s.   \label{ERK-y-2} 
\end{align}
Solving these recursions leads to
\begin{align}
    y_{n+1} &=  \tau \sum_{k=0}^{n} \mathrm{e}^{(n-k)\tau A}  \sum_{i=1}^{s} b_i (\tau A) g(t_{ki}),\label{ERK-y-recursion-1}  \\
    Y_{ni} &= \tau \sum_{k=0}^{n-1} \mathrm{e}^{(n+c_i-k-1)\tau A}  \sum_{m=1}^{s} b_m (\tau A) g(t_{km}) +\tau \sum_{j=1}^{i-1} a_{ij} (\tau A) g(t_{nj}),\quad i=1,\ldots,s. \label{ERK-y-recursion-2} 
\end{align}
By \eqref{solution-x} and \eqref{ERK-y-recursion-1}-\eqref{ERK-y-recursion-2}, the solution $u=x+y$ at $t=t_{n+1}$ and $t=t_{ni}$, respectively, can be written as
\begin{align}
u(t_{n+1}) & =x(t_{n+1}) + \tau \sum_{k=0}^{n} \mathrm{e}^{(n-k)\tau A}  \sum_{i=1}^{s} b_i (\tau A) g(t_{ki}) + y(t_{n+1}) - y_{n+1}, \label{Eqn:utn-recursion} \\
u(t_{ni}) & =x( t_{ni}) +  \tau \sum_{k=0}^{n-1} \mathrm{e}^{(n+c_i-k-1)\tau A}  \sum_{m=1}^{s} b_m (\tau A) g(t_{km}) +\tau \sum_{j=1}^{i-1} a_{ij} (\tau A) g(t_{nj}) \notag \\
&\quad  + y(t_{ni}) - Y_{ni}. \label{Eqn:utni-recursion} 
\end{align}

\begin{lemma}\label{Lem:estimate-yni}
Let the initial value problem \eqref{Eqn:semilinear-probolems} with $u_0\in X_\alpha$ satisfy Assumptions \ref{Ass:sectorial-A}-\ref{Ass:sectorial-Aalpha} and \ref{Ass:f-lip}. Consider for its numerical solution the EERK method \eqref{ERK-semilinear} satisfying \eqref{Eqn:linear-first-order}. Then the errors of \eqref{ERK-y-1} and \eqref{ERK-y-2} satisfy
$$
\begin{aligned}
\|y_n-y(t_n)\|_{\alpha}  & \leq C\tau^{1-\alpha}( \left|\log \tau\right| +1),\\
\|Y_{ni}-y(t_{ni})\|_{\alpha}  & \leq C\tau^{1-\alpha}( \left|\log \tau\right| +1) ,
\end{aligned}
$$
uniformly in $0\leq t_n  \leq t_{ni}\leq T$, respectively. The above constants $C$ depend on $T$, but are independent of $n$ and $\tau$.
\end{lemma}
\begin{proof}
By Lemma \ref{Lem:derivative} the time derivative of $g(t)=f(t,u(t))$ satisfies (almost everywhere)
$$
\|g^\prime(t)\| \leq \bigg\|\frac{\partial  f}{\partial t}(t,u(t))\bigg\| + \bigg\|\frac{\partial  f}{\partial u}(t,u(t))u^\prime(t)\bigg\| \leq Ct^{-1},\quad 0<  t \leq T
$$
for some $C$ depending on $T$. 
Since \eqref{Eqn:y} satisfies the assumptions in Theorem \ref{Thm:linear-first}, we have 
\begin{equation}\label{Eqn:yn-ytn}
\|y_n -y(t_n)\|_{\alpha} \leq  C\tau^{1-\alpha}( \left|\log \tau\right| +1).
\end{equation}
The mild solution of \eqref{Eqn:y} at $t=t_{ni}$ is given as
\begin{equation}\label{Eqn:exact-yni}
y(t_{ni})=\mathrm{e}^{c_i\tau A}y(t_{n}) + \int_0^{c_i\tau} \mathrm{e}^{(c_i\tau-\xi)A} g(t_n+\xi)\,\mathrm{d}\xi.
\end{equation}
Subtracting \eqref{Eqn:exact-yni} from \eqref{ERK-y-2} yields
$$
Y_{ni}-y(t_{ni})  =  \mathrm{e}^{c_i\tau A} (y_n -y(t_{n}) ) + \tau \sum_{j=1}^{s} a_{ij} (\tau A)g(t_{nj})  - \int_0^{c_i\tau} \mathrm{e}^{(c_i\tau-\xi)A} g(t_n+\xi)\,\mathrm{d}\xi.
$$
Estimating the above terms yields
$$
\begin{aligned}
\|Y_{ni}-y(t_{ni})\|_{\alpha} &  \leq C\|y_n -y(t_{n})\|_{\alpha} + C\tau^{1-\alpha} \sum_{j=1}^{s} \|\tau^{\alpha} a_{ij} (\tau A)\|_{L(X,D_A(\alpha,1))}\|g(t_{nj})\|  \\
&\quad + C\int_0^{c_i\tau} (c_i\tau-\xi)^{-\alpha}\|(c_i\tau-\xi)^{\alpha} \mathrm{e}^{(c_i\tau-\xi)A}\|_{L(X,D_A(\alpha,1))} \|g(t_n+\xi)\| \,\mathrm{d}\xi \\
& \leq C\|y_n -y(t_{n})\|_{\alpha} +C\tau^{1-\alpha},
\end{aligned}
$$
and the proof is completed by using \eqref{Eqn:yn-ytn}.  
\end{proof}
\begin{remark}
If the initial data satisfies $u_0\in D_A(\gamma,\infty)$ for some $\gamma \in(\alpha,1)$, then the error bounds in Lemma~\ref{Lem:estimate-yni} are independent of $\left| \log \tau \right|$.
\end{remark}

For first order methods, we have the following convergence result.

\begin{theorem}\label{Thm:semilinear-first}
For $\alpha \in [0,1)$, let the initial value problem \eqref{Eqn:semilinear-probolems} with $u_0\in X_\alpha$ satisfy Assumptions \ref{Ass:sectorial-A}-\ref{Ass:sectorial-Aalpha} and \ref{Ass:f-lip}. Consider for its numerical solution the EERK method \eqref{ERK-semilinear} satisfying \eqref{Eqn:linear-first-order}. Then, the error is bounded for $\tau$ sufficiently small by
\begin{equation*}
\|u_n -u(t_n)\|_{\alpha} \leq  C\tau^{1-\alpha}( \left|\log \tau\right| +1).
\end{equation*}
If further $u_0\in X_\gamma$ for some $\gamma \in (\alpha,1)$ and the part of $A$ in $X_\gamma$ is also sectorial in $X_\gamma$, then the error satisfies
$$
\begin{aligned}
&\|u_n -u(t_n)\|_{\alpha} \leq  C\tau^{1-\alpha} , \\
&\|u_n -u(t_n)\|_{\gamma} \leq  C\tau^{1-\gamma}.
\end{aligned}
$$
The above constants $C$ depend on $T$, but are independent of $n$ and $\tau$.
\end{theorem}
\begin{proof}
Solving the recursion \eqref{ERK-semilinear-1}-\eqref{ERK-semilinear-2} gives
\begin{align}
u_{n+1} & = \mathrm{e}^{(n+1)\tau A}u_0 + \tau \sum_{k=0}^{n}\mathrm{e}^{(n-k)\tau A}\sum_{i=1}^s b_i(\tau A)F_{ki}, \label{Eqn:un-recursion} \\
U_{ni} & = \mathrm{e}^{(n+c_i)\tau A}u_0
+ \tau \sum_{k=0}^{n-1}\mathrm{e}^{(n+c_i-k-1)\tau A}\sum_{m=1}^s b_m(\tau A)F_{km} + \tau \sum_{j=1}^{i-1} a_{ij}(\tau A)F_{nj}. \label{Eqn:Uni-recursion}
\end{align}
Subtracting \eqref{Eqn:utni-recursion} from \eqref{Eqn:Uni-recursion} yields
\begin{align}
%u_{n+1}-u(t_{n+1}) & =  \tau \sum_{k=0}^{n}\mathrm{e}^{(n-k)\tau L }\sum_{i=1}^s b_i(\tau L) ( F_{ki}  -g(t_{ki})) +(y_n-y(t_n))  , \label{Eqn:en-recursion} \\
U_{ni}-u(t_{ni}) & =   \tau \sum_{k=0}^{n-1}\mathrm{e}^{(n+c_i-k-1)\tau A }\sum_{m=1}^s b_m(\tau A)(F_{km}-g(t_{km}))+ \tau \sum_{j=1}^{i-1} a_{ij}(\tau A)(F_{nj}-g(t_{nj}))   \notag \\
&\quad 
+ Y_{ni}-y(t_{ni}) . \label{Eqn:eni-recursion} 
\end{align}
Let $E_{ni}=U_{ni}-u(t_{ni})$ denote the difference between the numerical and the exact solution and let $\|E_{n}\|_{\alpha} =\max_{1\leq j \leq s}\|E_{nj}\|_{\alpha}$. For $n=0$, we have
$$
\|E_{0i}\|_{\alpha} \leq C \tau^{1-\alpha} \max_{1\leq j \leq i-1}\|E_{0j}\|_{\alpha}+\|Y_{0i}-y(t_{0i})\|_{\alpha}.
$$
 With the help of Lemma \ref{Lem:estimate-yni}, we arrive at
 $$
 \|E_{0}\|_{\alpha} \leq C \tau^{1-\alpha}( \left|\log \tau\right| +1) .
 $$
For $n\geq 1$, estimating \eqref{Eqn:eni-recursion}  in $X_\alpha$ gives
\begin{align}
%\|u_n-u(t_n)\|_V & \leq  C\tau \sum_{j=0}^{n-2}t_{n-j-1}^{-\alpha}\sum_{i=1}^s \|U_{ji}-u(t_{j,i})\|_V + C\tau^{1-\alpha}\sum_{i=1}^s\|U_{n-1,i}-u(t_{n-1,i})\|_V  \notag \\
%&\quad + \|y_n-y(t_n)\|_V , \label{Eqn:en-recursion} \\
\|E_{ni}\|_{\alpha} & \leq  C\tau\sum_{k=0}^{n-2}t_{n-k-1}^{-\alpha}  \|E_{k}\|_{\alpha} + C\tau^{1-\alpha} \|E_{n-1}\|_{\alpha} 
+ C\tau^{1-\alpha} \max_{1 \leq j \leq i-1} \|E_{nj}\|_{\alpha}
+\|Y_{ni}-y(t_{ni})\|_{\alpha}. \notag  
\end{align}
Applying Lemma \ref{Lem:estimate-yni} again, we obtain
\begin{equation*}
\|E_n\|_{\alpha} \leq C\tau\sum_{k=0}^{n-2}t_{n-k-1}^{-\alpha}  \|E_{k}\|_{\alpha} + C\tau^{1-\alpha} \|E_{n-1}\|_{\alpha}  +C\tau^{1-\alpha}( \left|\log \tau\right| +1)  .
\end{equation*}
An application of a discrete Gronwall lemma \cite[Lemma 4]{HochbruckOstermann05:323} leads to 
\begin{equation}\label{Eqn:semilinear-stage-first}
\|E_n\|_{\alpha} \leq C\tau^{1-\alpha}( \left|\log \tau\right| +1),\quad 0\leq t_n < T.
\end{equation}
Note that $\|u_n -u(t_n)\|_{\alpha} \leq \|E_n\|_\alpha$ for $0\leq t_n < T$, which is a consequence of $c_1=0$. Further, subtracting \eqref{Eqn:utn-recursion} from \eqref{Eqn:un-recursion} gives
\begin{equation}\label{Eqn:en-recursion} 
u_{n+1}-u(t_{n+1})  =  \tau \sum_{k=0}^{n}\mathrm{e}^{(n-k)\tau A }\sum_{i=1}^s b_i(\tau A) ( F_{ki}  -g(t_{ki})) + y_n-y(t_n).
\end{equation}
Estimating this term in ${X_\alpha}$ and using \eqref{Eqn:semilinear-stage-first} yields
$$
\|u_n -u(t_n)\|_{\alpha} \leq  C\tau^{1-\alpha}( \left|\log \tau\right| +1),\quad 
0\leq t_n \leq T,
$$
which proves the first statement of the theorem.

If further $u_0\in X_\gamma$ for some $\gamma \in (\alpha,1)$ and the part of $A$ in $X_\gamma$ is also sectorial in $X_\gamma$, then the solution satisfies $u \in L^\infty( (0,T);X_\gamma)\cap C((0,T];X_\gamma)$. From Lemma \ref{Lem:derivative} we have $\|u^\prime (t)\|_{\alpha} < Ct^{\gamma - \alpha -1}$ and thus $\|g^\prime(t)\|\leq Ct^{\gamma - \alpha -1}$. It follows that
\begin{align}
\|y_n-y(t_n)\|_{\alpha}&\leq C\tau^{1-\alpha} ,&  \|y_n-y(t_n)\|_{\gamma} & \leq C\tau^{1-\gamma}, && 0\leq t_n\leq T, \label{Eqn:ynytn-alpha-gamma} \\
 \|Y_{ni}-y(t_{ni})\|_{\alpha}&\leq C\tau^{1-\alpha},&   \|Y_{ni}-y(t_{ni})\|_{\gamma}  &\leq C\tau^{1-\gamma} ,&& 0 \leq t_{ni}\leq T. \label{Eqn:yniytni-alpha-gamma} 
\end{align}
Following the steps of the preceding proof, we can derive 
\begin{align}
\|u_n -u(t_n)\|_{\alpha}  & \leq  C\tau^{1-\alpha} , \quad 0\leq t_n \leq T ,\label{Eqn:unutn-alpha-improved}  \\
\|U_{ni}-u(t_{ni})\|_{\alpha} & \leq  C\tau^{1-\alpha} , \quad 0\leq t_{ni} \leq T.\label{Eqn:uniutni-alpha-improved}  
\end{align}
Estimating \eqref{Eqn:en-recursion} in $X_\gamma$ gives
\begin{equation}\label{Eqn:extension-gamma}
\|e_{n}\|_{\gamma}   \leq  C\tau\sum_{k=0}^{n-2}t_{n-k-1}^{-\gamma}  \|E_{k}\|_{\alpha} + C\tau^{1-\gamma} \|E_{n-1}\|_{\alpha} 
+\|y_{n}-y(t_{n})\|_{\gamma}. 
\end{equation}
Plugging \eqref{Eqn:ynytn-alpha-gamma} and \eqref{Eqn:uniutni-alpha-improved}  into the above estimate yields
$$
\|u_n -u(t_n)\|_{\gamma} \leq  C\tau^{1-\gamma} ,\quad 
0\leq t_n \leq T,
$$
and the proof is completed.
\end{proof}

Next we consider the error estimates for second order methods. Such methods satisfy the order conditions (see \cite{HochbruckOstermann05:1069,LuanOstermann13:3431})
\begin{equation}\label{Eqn:semilinear-second-order}
 \sum_{i=1}^s b_i(z) = \varphi_1(z),\quad \sum_{i=1}^s b_i(z)c_i = \varphi_2( z),\quad\mbox{and}\quad
\sum_{j=1}^{i-1} a_{ij}(z) =c_i\varphi_1(c_iz)~\mbox{for}~2\leq i\leq s.
 \end{equation}

\begin{lemma}\label{Lem:estimate-yni-second}
For $\alpha\in[0,1)$, let the initial value problem \eqref{Eqn:semilinear-probolems} satisfy Assumptions \ref{Ass:sectorial-A}-\ref{Ass:sectorial-Aalpha} and \ref{Ass:f-lip}. Assume that $Au_0+f(0,u_0)\in X_\alpha$ and let $f:[0,T]\times X_\alpha \to X$ be Fr\'echet differentiable with a locally Lipschitz continuous derivative.
Consider for the numerical solution of \eqref{Eqn:semilinear-probolems} the EERK method \eqref{ERK-semilinear} satisfying \eqref{Eqn:semilinear-second-order}. Then, for any $\gamma \in [\alpha, 1)$ the errors of \eqref{ERK-y-1} and \eqref{ERK-y-2} satisfy
$$
\begin{aligned}
\|y_n-y(t_n)\|_{D_{A}(\gamma,1)}&  \leq C\tau^{2-\gamma}( \left|\log \tau\right| +1) , \\
\|Y_{ni}-y(t_{ni})\|_{D_{A}(\gamma,1)} & \leq C\tau^{2-\gamma}( \left|\log \tau\right| +1),
\end{aligned}
$$
uniformly in $0\leq t_n  \leq t_{ni}\leq T$, respectively. The above constants $C$ depend on $T$, but are independent of $n$ and $\tau$.
\end{lemma}
\begin{proof}
We consider the system
$$
\left\{\begin{aligned}
&u^\prime  =Au +  f(t,u),    && \,  u(0)=u_0, \\
&u^{\prime\prime}  =Au^\prime +  \frac{\partial f}{\partial t}(t,u) +\frac{\partial f}{\partial u}(t,u)u^\prime , &&\, u^\prime(0)=Au_0+f(0,u_0),
\end{aligned}\right.
$$
which can be reformulated as
\begin{equation}\label{Eqn:reformulation}
U^\prime=\mathcal{A}U + \mathcal{F}(t,U),\quad U_0 =[u_0,Au_0+f(0,u_0)]^T,
\end{equation}
where
$$
U=\begin{bmatrix} u \\ w \end{bmatrix} ,\quad
\mathcal{A}=  \begin{bmatrix}   A & 0 \\ 0 & A \end{bmatrix}  , \quad
\mathcal{F}(t,U)=  \begin{bmatrix}   f(t,u) \\ \frac{\partial f}{\partial t}(t,u) +\frac{\partial f}{\partial u}(t,u)w \end{bmatrix} .
$$
Note that $U_0 \in  D(A) \times X_\alpha\subset X_\alpha\times X_\alpha$ and $\mathcal{F}$ is locally Lipschitz bounded: for every $\widetilde{R}>0$, there exists $\widetilde{L}=\widetilde{L}(\widetilde{R},T)$ such that
\begin{equation*}
\|\mathcal{F}(t_1,V_1)-\mathcal{F}(t_2,V_2)\|_{X\times X} \leq \widetilde{L} \big( |t_1-t_2| + \|V_1-V_2\|_{X_\alpha\times X_\alpha} \big)
\end{equation*}
for all $t\in[0,T]$ and $V_1,V_2\in X_\alpha\times X_\alpha$ with $\max\big( \|V_1\|_{X_\alpha\times X_\alpha},\|V_2\|_{X_\alpha\times X_\alpha})\leq \widetilde{R}$. Therefore, there exists a unique mild solution of \eqref{Eqn:reformulation} satisfying
$$
\|U(t)\|_{X_\alpha\times X_\alpha} \leq C\|U_0\|_{X_\alpha\times X_\alpha} + CT^{1-\alpha}  \|\mathcal{F}(\cdot,U(\cdot))\|_{L^\infty(0,T;X\times X)}, \quad 0\leq t \leq T ,
$$
and 
$$
\|u^{\prime\prime}(t)\|_{\alpha} \leq \|U^\prime(t)\|_{X_\alpha\times X_\alpha} \leq Ct^{-1},\quad 0<t\leq T.
$$
Then the second order time derivative of $g(t)=f(t,u(t))$ satisfies (almost everywhere)
$$
\|g^{\prime\prime}(t)\| \leq  Ct^{-1},\quad 0<  t \leq T.
$$
In a similar way in the proof of Lemma~\ref{Lem:estimate-yni}, the desired order is thus obtained with the help of Theorem~\ref{Thm:linear-second}.
\end{proof}

With the help of Lemma~\ref{Lem:estimate-yni-second}, the following convergence result for second order methods can be obtained.
\begin{theorem}\label{Thm:semilinear-second}
For $\alpha\in [0,1)$, let the initial value problem \eqref{Eqn:semilinear-probolems} satisfy Assumptions \ref{Ass:sectorial-A}-\ref{Ass:sectorial-Aalpha} and \ref{Ass:f-lip}. Assume that $Au_0+f(0,u_0)\in X_\alpha$ and let $f:[0,T]\times X_\alpha\to X$ be Fr\'echet differentiable with a locally Lipschitz continuous derivative.
Consider for the numerical solution of \eqref{Eqn:semilinear-probolems} the EERK method \eqref{ERK-semilinear} satisfying \eqref{Eqn:semilinear-second-order}. Then, for any $\gamma\in[\alpha,1)$ and $\tau$ sufficiently small the error is bounded by
%$$
%\begin{aligned}
%&\|u_n -u(t_n)\|_{\alpha} \leq  C\tau^{2-\alpha}( \left|\log \tau\right| +1), \\
%&\|u_n -u(t_n)\|_{D_A(\gamma,1)} \leq  C\tau^{2-\gamma}( \left|\log \tau\right| +1),\quad 
%\gamma\in[\alpha,1).
%\end{aligned}
%$$
$$
\|u_n -u(t_n)\|_{D_A(\gamma,1)} \leq  C\tau^{2-\gamma}( \left|\log \tau\right| +1),
$$
where the constant $C$ depends on $T$, but is independent of $n$ and $\tau$.
\end{theorem}

\section{Examples}\label{Sec:example}
In this section, we shall use the nonsmooth data error estimates to study some numerical solutions of semilinear parabolic equations.

\begin{example}\label{Exa:allen-cahn}
We consider the following Allen--Cahn equation:  
\begin{equation}\label{Eqn:allen-cahn}
\left\{\begin{aligned}
&\partial_t u =\varepsilon^2 \Delta u + u-u^3,   && \,  t>0, \ x \in  \Omega,\\
&u(0,x)=u_0(x),   && \, x \in \Omega, \\
&u(t,x)=0,  && \, t>0,\ x \in   \partial\Omega,
\end{aligned}\right.
\end{equation}
where $\varepsilon>0$ and $\Omega$ is an open bounded set in $\mathbb{R}^n$ with sufficiently smooth boundary $\partial\Omega$.
\end{example}

We choose $X=C(\overline{\Omega})$ and $\alpha=0$. The realization $A:D(A)\to X$ of the operator $\varepsilon^2 \Delta$ in $X$, with domain 
$$
D(A)=\{u\in C(\overline{\Omega}):  \Delta u \in C(\overline{\Omega}),~u|_{\partial\Omega}=0\},
$$
is sectorial in $X$ due to \cite[Corollary 3.1.21]{Lunardi95Book}. The nonlinearity $f(u)=u-u^3$ is locally Lipschitz continuous
$$
\|f(v_1)-f(v_2)\|_{C(\overline{\Omega})} \leq L(R)\|v_1-v_2\|_{C(\overline{\Omega})},\quad  v_1,v_2\in C(\overline{\Omega}),
$$
where $L(R)=1 +3R^2$ and $R=\max\big\{ \|v_1\|_{C(\overline{\Omega})},\|v_2\|_{C(\overline{\Omega})}\big\}$. If  $u_0\in C(\overline{\Omega})$, this problem admits a unique mild solution for some $T>0$. Consider the numerical solution of \eqref{Eqn:allen-cahn} by the EERK method \eqref{ERK-semilinear} satisfying \eqref{Eqn:linear-first-order}. Then, the error for $\tau$ sufficiently small is bounded  by
$$
\|u_n-u(t_n)\|_{C(\overline{\Omega})} \leq C\tau(\left|\log \tau \right| +1).
$$

Now we show that Lemma \ref{Lem:estimate-yni} can be improved for this example. From \cite[Theorem 3.1.29]{Lunardi95Book}, it holds
\begin{equation}\label{Eqn:DAinfty}
D_A(\gamma,\infty) = C^{2\gamma}_0(\overline{\Omega}),\quad \gamma \in (0,\tfrac12)\cup (\tfrac12,1),
\end{equation}
and from \cite[Proposition 3.1.27 and Proposition 3.1.28]{Lunardi95Book}, we get
\begin{equation*}%\label{Eqn:DAinfty-12}
D_A(\tfrac12,1) \subset  C^{1}_0(\overline{\Omega}) \subset D_A(\tfrac12,\infty).
\end{equation*}
Thanks to \cite[Proposition 2.2.7 and Theorem 3.1.25]{Lunardi95Book}, the part of $A$ in $C_0^{2\gamma}(\overline{\Omega})$ is sectorial in $C_0^{2\gamma}(\overline{\Omega})$ for $0<\gamma <1$ . 
If further $u_0 \in C_0^{2\gamma}(\overline{\Omega})$ for some $\gamma\in (0,1)$, we have
$$
\|u(t)\|_{C_0^{2\gamma}(\overline{\Omega})} \leq C,\quad 0\leq t \leq T,\qquad
\|u^\prime(t)\|_{C_0^{2\gamma}(\overline{\Omega})} < Ct^{-1},\quad 0\leq t \leq T,
$$
where the latter follows from Lemma \ref{Lem:derivative}.
Since $f:C_0^{2\gamma}(\overline{\Omega})\to C_0^{2\gamma}(\overline{\Omega})$ is continuous in this example, the function $g(t)=f(u(t))$ satisfies
$$
\|g(t)\|_{C_0^{2\gamma}(\overline{\Omega})} \leq C,\quad 0\leq t \leq T,\qquad
\|g^\prime(t)\|_{C_0^{2\gamma}(\overline{\Omega})} < Ct^{-1},\quad 0\leq t \leq T.
$$
Therefore, the results in Lemma \ref{Lem:estimate-yni} can be improved, which leads to the refined estimate
$$
\|u_n-u(t_n)\|_{C_0^{2\gamma}(\overline{\Omega})} \leq C\tau (\left|\log \tau \right| +1).
$$

Consider now the numerical solution of \eqref{Eqn:allen-cahn} by the EERK method \eqref{ERK-semilinear} satisfying \eqref{Eqn:linear-second-order}. 
If further $Au_0\in C(\overline{\Omega})$ (i.e., $u_0 \in D(A)$), then for any $\gamma \in [0,1)$ the error satisfies
$$
\|u_n-u(t_n)\|_{C_0^{2\gamma}(\overline{\Omega})} \leq C\tau^{2-\gamma}(\left|\log \tau \right| +1).
$$

\begin{example}\label{Exa:Burgers}
We consider the following Burgers’ equation:  
\begin{equation}\label{Eqn:Burgers}
\left\{\begin{aligned}
&\partial_t u = \nu \Delta u - u \partial_x u - u \partial_y u,   &&  \, t>0,\ x \in  \Omega,\\
&u(0,x)=u_0(x), \!\!\!  && \, x \in \Omega, \\
&u(t,x)=0, \!\!\!  && \, t>0,\ x \in   \partial\Omega,
\end{aligned}\right.
\end{equation}
where $\nu > 0$ is the viscosity and $\Omega$ is an open bounded set in $\mathbb{R}^n$ with sufficiently smooth boundary $\partial\Omega$.
\end{example}

We still choose $X=C(\overline{\Omega})$. The operator $A$ can be defined in the same manner as before, except that $\varepsilon^2$ is replaced with $\nu$. The domains of the operator remain the same. 
Let $f(u)=- u \partial_x u - u \partial_y u$. The maximal domain of $f$ in $C(\overline{\Omega})$ is $X_{1/2}=C_0^1(\overline{\Omega})$. Note that
\begin{equation}\label{Eqn:exa2-f}
\|f(v_1)-f(v_2)\|_{C(\overline{\Omega})} \leq L\|v_1-v_2\|_{C_0^1(\overline{\Omega})},\quad v_1,v_2\in C_0^1(\overline{\Omega}),
\end{equation}
where $L(R)=2R$ and $R=\max\big\{ \|v_1\|_{C_0^1(\overline{\Omega})},\|v_2\|_{C_0^1(\overline{\Omega})}\big\}$. If  $u_0\in C_0^1(\overline{\Omega})$, this problem admits a unique mild solution with some $T>0$. Consider the numerical solution of \eqref{Eqn:Burgers} by the EERK method \eqref{ERK-semilinear} satisfying \eqref{Eqn:linear-first-order}. If $u_0\in C_0^{2\gamma}(\overline{\Omega})$ for some $\gamma \in[\tfrac12,1)$, then the error for $\tau$ sufficiently small is bounded  by
$$
\|u_n-u(t_n)\|_{C_0^{2\gamma}(\overline{\Omega})} \leq C\tau^{1/2}(\left|\log \tau \right| +1).
$$
This estimate provides a sharper bound compared to Theorem \ref{Thm:semilinear-first}, by exploiting the continuity of $f:C_0^{2\gamma}(\overline{\Omega})\to C_0^{2\gamma-1}(\overline{\Omega})$.

If further the compatibility condition $Au_0+f(u_0) \in C_0^1(\overline{\Omega})$ holds, then for any $\gamma\in [\tfrac12,1)$ the error of the numerical solution of \eqref{Eqn:Burgers} by the EERK method \eqref{ERK-semilinear} satisfying \eqref{Eqn:linear-second-order} is bounded by
$$
\|u_n-u(t_n)\|_{C_0^{2\gamma}(\overline{\Omega})} \leq C\tau^{2-\gamma}(\left|\log \tau \right| +1).
$$

\section{Application: split exponential integrators}\label{Sec:split}
As an application of the nonsmooth data error estimates, we consider split exponential integrators for the problem \eqref{Eqn:semilinear-probolems}. For this application, let $A=A_1+A_2$, where $A_1$ and $A_2$ are assumed to satisfy the following assumption.

\begin{assumption}\label{Ass:sectorial-AB}
Let $A_1: D(A_1) \rightarrow X$ and $A_2: D(A_2) \rightarrow X$ be sectorial and let the corresponding generated analytic semigroups $\{\mathrm{e}^{tA_1}\}$ and $\{\mathrm{e}^{tA_2}\}$ commute in the sense of $\mathrm{e}^{tA_1}\mathrm{e}^{tA_2}=\mathrm{e}^{tA_2}\mathrm{e}^{tA_1}$.
\end{assumption}
Under this assumption, it holds $\mathrm{e}^{t(A_1+A_2)}=\mathrm{e}^{tA_1}\mathrm{e}^{tA_2}$ and $D(A)=\overline{D(A_1)\cap D(A_2)}$.
%This assumption infers the parabolic smoothing property
%\begin{equation}\label{Est:etA1etA2}
%\|\mathrm{e}^{tA_1}\|  + \|\mathrm{e}^{tA_2}\|  \leq C,\quad t \geq 0.
%\end{equation}
Applying the second order EERK in \cite[Equation (5.3)]{HochbruckOstermann05:1069} to the problem \eqref{Eqn:semilinear-probolems} gives 
\begin{subequations}\label{ERK-original} 
\begin{align}
    \widetilde{u}_{n+1} &= \mathrm{e}^{\tau A} \widetilde{u}_n +   \tau \widetilde{b}_1(\tau A)  \widetilde{F}_{n1} + \tau \widetilde{b}_2(\tau A) \widetilde{F}_{n2} , \label{ERK-original-1} \\
    \widetilde{U}_{n1}&=\widetilde{u}_{n},\quad \widetilde{U}_{n2} = \mathrm{e}^{c_2\tau A} \widetilde{u}_n +\tau \widetilde{a}_{21}(\tau A)  \widetilde{F}_{n1} ,  \label{ERK-original-2} \\
    \widetilde{F}_{nj} &= f(t_{nj}, \widetilde{U}_{nj}),\quad j=1,2,  \label{ERK-original-3} 
\end{align}
\end{subequations}
where the coefficients are given by
$$
\begin{aligned}
\widetilde{b}_1(\tau A) & = \varphi_1(\tau A)-\tfrac{1}{c_2}\varphi_2(\tau A) ,  \\
\widetilde{b}_2(\tau A) & =\tfrac{1}{c_2}\varphi_2(\tau A) , \\
\widetilde{a}_{21}(\tau A ) & = c_2 \varphi_1(c_2\tau A).
\end{aligned}
$$

Recently, the following split version of integrator \eqref{ERK-original}, referred to as ERK2L, was constructed in \cite{CaliariCassini25},
\begin{subequations}\label{ERK-split} 
\begin{align}
    {u}_{n+1} &= \mathrm{e}^{\tau (A_1+A_2)} {u}_n + \tau b_1(\tau A_1,\tau A_2) {F}_{n1} + \tau b_2(\tau A_1,\tau A_2) {F}_{n2} , \label{ERK-split-1} \\
    {U}_{n1} &=u_n,\quad {U}_{n2} = \mathrm{e}^{c_2\tau  (A_1+ A_2)} {u}_n + \tau a_{21}(\tau A_1,\tau A_2) {F}_{n1} ,  \label{ERK-split-2} \\
    {F}_{nj} &= f(t_{nj}, {U}_{nj}),\quad j=1,2,  \label{ERK-split-3} 
\end{align}
\end{subequations}
where the coefficients are given by
$$
\begin{aligned}
b_1(\tau A_1,\tau A_2) & = \varphi_1(\tau  A_1)  \varphi_1(\tau A_2) -\tfrac{2}{c_2}\varphi_2(\tau A_1) \varphi_2(\tau A_2),  \\
b_2(\tau A_1,\tau A_2) & =\tfrac{2}{c_2}\varphi_2(\tau A_1) \varphi_2(\tau A_2), \\
a_{21}(\tau A_1,\tau A_2) & = c_2\varphi_1(c_2\tau A_1)\varphi_1(c_2\tau A_2).
\end{aligned}
$$
The choice of these coefficients is motivated by
\begin{equation}\label{Eqn:split-approximations}
\varphi_1(\tau(A_1+A_2)) \approx \varphi_1(\tau A_1) \varphi_1(\tau A_2),
\quad \varphi_2(\tau(A_1+A_2)) \approx 2\varphi_2(\tau A_1) \varphi_2(\tau A_2).
\end{equation}
For the case $\alpha=0$, it has been proven in \cite{CaliariCassini25} that under strong smoothness assumptions on the solution and the nonlinearity $f$, the error of ERK2L is bounded by 
$$
\|e_n\| \leq \left\{\begin{aligned} 
&C\tau^2(\left|\log \tau\right|+1), && \mbox{if}~f(t,u(t))\in D(A_1)~\mbox{or}~f(t,u(t))\in D(A_2),\\
& C\tau, && \mbox{else}. \end{aligned}\right.
$$
In this section, we conduct a more refined convergence analysis. As a first result we show that the approximation errors of \eqref{Eqn:split-approximations} satisfy the following bound. 

\begin{lemma}\label{Est:split-varphi}
Under Assumptions \ref{Ass:sectorial-A} and \ref{Ass:sectorial-AB}, it holds, for $k\geq 1$ and $0<t\leq T$,
\begin{itemize}
\item[\rm (a)] $\big\| t^{-\beta_1+\beta_2} \big( \varphi_k(t (A_1+A_2))  - k! \varphi_k(t A_1)\varphi_k(t A_2) \big) \big\|_{L(D_{A_1}(\beta_1,\infty),D_A(\beta_2,1))}  \leq C,\quad 0\leq \beta_1,\beta_2 <1$, 
\item[\rm (b)] $\big\| t^{-1+\beta_2} \big( \varphi_k(t(A_1+A_2))  - k! \varphi_k(t A_1)\varphi_k(t A_2) \big) \big \|_{L(D(A_1),D_A(\beta_2,1))}  \leq C,\quad 0\leq  \beta_2 <1$,
\end{itemize}
where the constant $C$ depends on $k$ and $T$, but is independent of $t$.
\end{lemma}
\begin{proof}
By using the relation 
$$
 \mathrm{e}^{t(1-\xi)A_1}= \mathrm{e}^{t(1-\eta)A_1} -\int_\eta^\xi t A_1\mathrm{e}^{t(1-\sigma)A_1} \,\mathrm{d}\sigma,
$$
we obtain
\begin{align}
&\varphi_k(t(A_1+A_2)) - k!\varphi_k(t A_1)\varphi_k(t A_2) \notag  \\
&\qquad = \int_0^1 \mathrm{e}^{t(1-\eta)(A_1+A_2)}\frac{\eta^{k-1}}{(k-1)!}\,\mathrm{d}\eta -\int_0^1 \int_0^1\mathrm{e}^{t(1-\eta)A_2}\mathrm{e}^{t(1-\xi)A_1}\frac{k(\eta\xi)^{k-1}}{(k-1)!}\, \mathrm{d}\eta  \,\mathrm{d}\xi \notag \\
&\qquad=  \frac{k}{(k-1)!}\int_0^1\int_0^1 \int_\eta^\xi t  \mathrm{e}^{t(1-\eta)A_2} A_1\mathrm{e}^{t(1-\sigma)A_1} (\eta\xi)^{k-1} \,\mathrm{d}\sigma \, \mathrm{d}\eta  \,\mathrm{d}\xi \notag   \\
&\qquad=  -\frac{k}{(k-1)!}\int_0^1\int_{\xi}^1 \int_\xi^\eta t \mathrm{e}^{t(1-\eta)(A_1+A_2)}  A_1 \mathrm{e}^{t(\eta-\sigma)A_1} (\eta\xi)^{k-1} \,\mathrm{d}\sigma \, \mathrm{d}\eta  \,\mathrm{d}\xi \notag \\
&\quad\qquad  +\frac{k}{(k-1)!}\int_0^1\int_0^\xi \int_\eta^\xi t \mathrm{e}^{t(1-\xi)(A_1+A_2)}\mathrm{e}^{t(\xi-\eta)A_2}A_1\mathrm{e}^{t(\xi-\sigma)A_1}  (\eta\xi)^{k-1} \,\mathrm{d}\sigma \, \mathrm{d}\eta  \,\mathrm{d}\xi \notag.
\end{align}
Recalling \eqref{est:X-DAalpha1}, property (b) of Lemma \ref{Lem:est-etA} and the boundedness of the semigroup yields
\begin{align}
&\|\varphi_k(t(A_1+A_2)) - k!\varphi_k(t A_1)\varphi_k(t A_2)\|_{L(D_{A_1}(\beta_1,\infty),D_A(\beta_2,1))} \notag  \\
&\qquad \leq  C\int_0^1\int_{\xi}^1 \int_\xi^\eta t \|\mathrm{e}^{t(1-\eta)(A_1+A_2)}\|_{L(X,D_A(\beta_2,1))}  \|A_1 \mathrm{e}^{t(\eta-\sigma)A_1}\|_{L(D_{A_1}(\beta_1,\infty),X)}   \,\mathrm{d}\sigma \, \mathrm{d}\eta  \,\mathrm{d}\xi \notag \\
&\quad\qquad  +C\int_0^1\int_0^\xi \int_\eta^\xi t\|\mathrm{e}^{t(1-\xi)(A_1+A_2)}\|_{L(X,D_A(\beta_2,1))}  \|A_1\mathrm{e}^{t(\xi-\sigma)A_1}\|_{L(D_{A_1}(\beta_1,\infty),X)}   \,\mathrm{d}\sigma \, \mathrm{d}\eta  \,\mathrm{d}\xi \notag \\
&\qquad \leq Ct^{\beta_1-\beta_2} \notag
\end{align}
and
\begin{align}
&\|\varphi_k(t(A_1+A_2)) - k!\varphi_k(t A_1)\varphi_k(t A_2)\|_{L(D(A_1),D_A(\beta_2,1))} \notag  \\
&\qquad \leq  C\int_0^1\int_{\xi}^1 \int_\xi^\eta t  \|\mathrm{e}^{t(1-\eta)(A_1+A_2)}\|_{L(X,D_A(\beta_2,1))}  \|A_1\|_{L(D(A_1),X)}\,\mathrm{d}\sigma \, \mathrm{d}\eta  \,\mathrm{d}\xi \notag \\
&\quad\qquad  +C\int_0^1\int_0^\xi \int_\eta^\xi t\|\mathrm{e}^{t(1-\xi)(A_1+A_2)}\|_{L(X,D_A(\beta_2,1))} \|A_1\|_{L(D(A_1),X)}   \,\mathrm{d}\sigma \, \mathrm{d}\eta  \,\mathrm{d}\xi \notag \\
&\qquad \leq Ct^{1-\beta_2}, \notag
\end{align}
which concludes the proof.
\end{proof}
\begin{remark}
The following estimates
\begin{itemize}
\item[\rm (a)] $\big\| t^{-\beta_1+\beta_2} \big( \varphi_k(t (A_1+A_2))  - k! \varphi_k(t A_1)\varphi_k(t A_2) \big) \big\|_{L(D_{A_2}(\beta_1,\infty),D_A(\beta_2,1))}  \leq C,\quad 0\leq \beta_1,\beta_2 <1$, 
\item[\rm (b)] $\big\| t^{-1+\beta_2} \big( \varphi_k(t(A_1+A_2))  - k! \varphi_k(t A_1)\varphi_k(t A_2) \big) \big \|_{L(D(A_2),D_A(\beta_2,1))}  \leq C,\quad 0\leq  \beta_2 <1$,
\end{itemize}
also hold true due to the fact that $\mathrm{e}^{tA_1}\mathrm{e}^{tA_2}=\mathrm{e}^{tA_2}\mathrm{e}^{tA_1}$.
\end{remark}
Using property (a) in Lemma \ref{Est:split-varphi} (with $\beta_1=0$, $\beta_2=\beta$) and \eqref{Est:varphik-X-DAalpha1}, we obtain
\begin{align}\label{Est:varphikvarphik-X-DAalpha1}
 \| t^\beta   \varphi_k( t A_1)\varphi_k( t A_2) \|_{L(X, D_A(\beta, 1))} 
&  \leq  \big\| t^{\beta}\big( (k!)^{-1} \varphi_k(t (A_1+A_2))  -  \varphi_k(t A_1)\varphi_k(t A_2) \big) \big\|_{L(X,D_A(\beta,1))} 
 \notag \\
 & \quad+   \big\| t^{\beta}(k!)^{-1}  \varphi_k(t (A_1+A_2)) \big\|_{L(X,D_A(\beta,1))}  \notag \\
 & \leq  C,\quad 0\leq  \beta < 1,~0<t \leq T.
\end{align}
In addition, we need the following relation. Its proof can be found in \cite{CaliariCassini25}.
\begin{lemma}\label{Lem:split-difference-representation}
Under Assumptions \ref{Ass:sectorial-A} and \ref{Ass:sectorial-AB}, it holds, for $k\geq 1$,
$$
\begin{aligned}
&\varphi_k(t(A_1+A_2)) - k!\varphi_k(t A_1)\varphi_k(t A_2) \\
&\qquad  =   \frac{k}{(k-1)!}\int_0^1\int_0^1 \int_\eta^\xi t^2(1-\eta)  A_1\mathrm{e}^{t(1-\sigma)A_1} A_2\varphi_1(t(1-\eta)A_2)(\eta\xi)^{k-1} \,\mathrm{d}\sigma \, \mathrm{d}\eta  \,\mathrm{d}\xi .
\end{aligned}
$$
\end{lemma}

We first consider the nonsmooth error estimate for split version of the exponential Euler method: 
\begin{equation}\label{Eqn:split-euler}
u_{n+1}=\mathrm{e}^{\tau(A_1+A_2)}u_n + \tau \varphi_1(\tau A_1)\varphi_1(\tau A_2)f(t_n,u_n).
\end{equation}
\begin{theorem}\label{Thm:semilinear-split-first}
For $\alpha \in [0,1)$, let the initial value problem \eqref{Eqn:semilinear-probolems} with $u_0\in X_\alpha$ satisfy Assumptions \ref{Ass:sectorial-A}-\ref{Ass:sectorial-Aalpha}, \ref{Ass:f-lip} and \ref{Ass:sectorial-AB}. Consider for its numerical solution the exponential Euler method \eqref{Eqn:split-euler}. Then, the error for $\tau$ sufficiently small is bounded by
\begin{equation*}
\|u_n -u(t_n)\|_{\alpha} \leq  C\tau^{1-\alpha}( \left|\log \tau\right| +1).
\end{equation*}
If further $u_0\in X_\gamma$ for some $\gamma \in (\alpha,1)$ and the part of $A$ in $X_\gamma$ is also sectorial in $X_\gamma$, the error is bounded by
$$
\begin{aligned}
\|u_n -u(t_n)\|_{\alpha} &\leq  C\tau^{1-\alpha}, \\
\|u_n -u(t_n)\|_{\gamma}& \leq  C\tau^{1-\gamma}.
\end{aligned}
$$
The above constants $C$ depend on $T$, but are independent of $n$ and $\tau$. 
\end{theorem}
\begin{proof}
Based on \eqref{Eqn:utn-recursion}, we rewrite the solution as
\begin{align}
u(t_{n+1})  &=x(t_{n+1}) + \tau \sum_{k=0}^{n} \mathrm{e}^{(n-k)\tau A}  \varphi_1 (\tau A_1,\tau A_2) g(t_{k}) + y(t_{n+1}) - y_{n+1}\notag  \\
&\quad  + \tau \sum_{k=0}^{n} \mathrm{e}^{(n-k)\tau A} \big( \varphi_1 (\tau (A_1+A_2)) - \varphi_1 (\tau A_1,\tau A_2) \big)g(t_{k}) .    \label{Eqn:split-euler-1}
\end{align}
Solving the recursion \eqref{Eqn:split-euler} gives
\begin{equation}\label{Eqn:split-euler-2}
u_{n+1}=\mathrm{e}^{(n+1)\tau A}+ \tau \sum_{k=0}^{n} \mathrm{e}^{(n-k)\tau A}  \varphi_1 (\tau A_1,\tau A_2) f(t_k,u_k).
\end{equation}
Subtracting \eqref{Eqn:split-euler-1} from \eqref{Eqn:split-euler-2} yields
\begin{align}
e(t_{n+1})  &=  \tau \sum_{k=0}^{n} \mathrm{e}^{(n-k)\tau A}  \varphi_1 (\tau A_1,\tau A_2) \big(f(t_k,u_k)-g(t_{k}) \big) - y(t_{n+1})+ y_{n+1} \notag  \\
&\quad  - \tau \sum_{k=0}^{n} \mathrm{e}^{(n-k)\tau A} \big( \varphi_1 (\tau (A_1+A_2)) - \varphi_1 (\tau A_1,\tau A_2) \big)g(t_{k})  .    \label{Eqn:split-euler-3}
\end{align}
By using the commutativity of the semigroups and Lemmas \ref{Lem:split-difference-representation} and \ref{Est:split-varphi}, we obtain the estimate of the term involving split approximations:
$$
\begin{aligned}
&\Bigg\|\tau \sum_{k=0}^{n} \mathrm{e}^{(n-k)\tau A} \big( \varphi_1 (\tau (A_1+A_2)) - \varphi_1 (\tau A_1,\tau A_2) \big)g(t_{k})\Bigg\|_{\alpha}  \\
&\qquad \leq \Bigg\|\tau^3 \sum_{k=1}^{n-3} \mathrm{e}^{(n-k)\tau A} \int_0^1\int_0^1 \int_\eta^\xi  (1-\eta)  A_1\mathrm{e}^{\tau(1-\sigma)A_1} A_2\varphi_1(\tau(1-\eta) A_2)  \,\mathrm{d}\sigma \, \mathrm{d}\eta  \,\mathrm{d}\xi \cdot g(t_{k})\Bigg\|_{\alpha}  \\
&\qquad \quad +\tau \sum_{k=n-2}^{n}  \big\|\mathrm{e}^{(n-k)\tau A} \big( \varphi_1 (\tau (A_1+A_2)) - \varphi_1 (\tau A_1,\tau A_2) \big)g(t_{k}) \big\|_{\alpha} \\
&\qquad \leq  C\tau^3 \sum_{k=0}^{n-3} \|\mathrm{e}^{\frac{1}{2} (n-k)\tau A}\|_{L(X,D_A(\alpha,1))}  \|A_1\mathrm{e}^{\frac{1}{2} (n-k)\tau A_1}\|\|A_2\mathrm{e}^{\frac{1}{2} (n-k)\tau A_2}\|  \|g(t_{k})\| \\
&\qquad\quad + C\tau \|\varphi_1 (\tau A) - \varphi_1 (\tau A_1,\tau A_2)\|_{L(X,D_A(\alpha,1))}\big(\|g(t_{n-2})\|+\|g(t_{n-1})\|+\|g(t_n)\|\big) \\
&\qquad\leq C\tau^3 \sum_{k=0}^{n-3} t_{n-k}^{-2-\alpha} + C\tau^{1-\alpha}\\
&\qquad\leq C\tau^{1-\alpha} .
\end{aligned}
$$
Recalling the estimate of $\|y_{n+1} - y(t_{n+1})\|_{\alpha}$ given in Lemma \ref{Lem:estimate-yni}, from \eqref{Eqn:split-euler-3} we have
$$
\|u_n -u(t_n)\|_{\alpha} \leq  C\tau^{1-\alpha}( \left|\log \tau\right| +1).
$$

If in addition $u_0\in X_\gamma$ for some $\gamma \in (\alpha,1)$ and the part of $A$ in $X_\gamma$ is also sectorial in $X_\gamma$, we can obtain the desired results in a similar manner as that presented in the proof of Theorem~\ref{Thm:semilinear-first}. The proof is completed.
\end{proof}

The convergence of ERK2L is presented in the following theorem. For convenience, we denote $D_{A_i}(1,\infty)=D(A_i)$.

\begin{theorem}\label{Thm:semilinear-split-second}
For $\alpha\in[0,1)$, let the initial value problem \eqref{Eqn:semilinear-probolems} satisfy Assumptions \ref{Ass:sectorial-A}-\ref{Ass:sectorial-Aalpha}, \ref{Ass:f-lip} and \ref{Ass:sectorial-AB}. Assume that $Au_0+f(0,u_0)\in X_\alpha$ and let $f:[0,T]\times X_\alpha \to X$ be Fr\'echet differentiable with a locally Lipschitz continuous derivative. 
Consider for the numerical solution the split integrator ERK2L \eqref{ERK-split}. If $f:[0,T] \times D({A}) \to D_{A_i}(\kappa,\infty)$ is continuous for some $\kappa \in [0,1]$ and $i=1$ or $2$, then for any $\gamma\in[\alpha,1)$ the error for $\tau$ sufficiently small is bounded by
\begin{equation*}\label{Eqn:semilinear-split-second}
\|u_n -u(t_n)\|_{D_A(\gamma,1)} \leq  C\tau^{1+\kappa-\gamma}(\left| \log \tau \right| +1),
\end{equation*}
where the constant $C$ depends on $T$, but is independent of $n$ and $\tau$.
\end{theorem}

\begin{proof}
Based on \eqref{Eqn:utn-recursion}-\eqref{Eqn:utni-recursion}, we rewrite the solution as
\begin{align}
u(t_{n+1}) & =x(t_{n+1}) + \tau \sum_{k=0}^{n} \mathrm{e}^{(n-k)\tau A}  \sum_{i=1}^{2} b_i (\tau A_1,\tau A_2) g(t_{ki}) +  y(t_{n+1}) - y_{n+1}   \notag \\
&\quad +\tau \sum_{k=0}^{n} \mathrm{e}^{(n-k)\tau A}  \sum_{i=1}^{2} \big(\widetilde{b}_i (\tau (A_1+A_2))-b_i (\tau A_1,\tau A_2) \big) g(t_{ki})  , \label{Eqn:split-utn-recursion} \\
u(t_{n2}) & =x( t_{n2}) +  \tau \sum_{k=0}^{n-1} \mathrm{e}^{(n+c_2-k-1)\tau A}  \sum_{m=1}^{2} b_m (\tau A_1,\tau A_2) g(t_{km}) +\tau   a_{21} (\tau A_1,\tau A_2) g(t_{n1}) \notag \\
&\quad  +  y(t_{n2}) - Y_{n2} +\tau \sum_{k=0}^{n-1} \mathrm{e}^{(n+c_2-k-1)\tau A}  \sum_{m=1}^{2} \big(\widetilde{b}_m (\tau (A_1+A_2) )-b_m (\tau A_1,\tau A_2) \big) g(t_{km}) \notag \\
&\quad  +\tau  \big( \widetilde{a}_{21} (\tau (A_1+A_2)) -a_{21} (\tau A_1,\tau A_2)  \big) g(t_{n1}). \label{Eqn:split-utni-recursion} 
\end{align}
Solving the recursion \eqref{ERK-split-1}-\eqref{ERK-split-2} gives
\begin{align}
u_{n+1} & = \mathrm{e}^{(n+1)\tau A}u_0 + \tau \sum_{k=0}^{n}\mathrm{e}^{(n-k)\tau A}\sum_{i=1}^2 b_i(\tau A_1,\tau A_2)F_{ki}, \label{Eqn:split-un-recursion} \\
U_{n2} & = \mathrm{e}^{(n+c_i)\tau A}u_0
+ \tau \sum_{k=0}^{n-1}\mathrm{e}^{(n+c_2-k-1)\tau A}\sum_{m=1}^2 b_m(\tau A_1,\tau A_2)F_{km} + \tau  a_{21}(\tau A_1,\tau A_2)F_{n1}. \label{Eqn:split-Uni-recursion}
\end{align}
Subtracting \eqref{Eqn:split-utn-recursion} from \eqref{Eqn:split-un-recursion} and \eqref{Eqn:split-utni-recursion} from \eqref{Eqn:split-Uni-recursion} yields
\begin{align}
e_{n+1}   & =  \tau \sum_{k=0}^{n}\mathrm{e}^{(n-k)\tau A }\sum_{i=1}^2 {b}_i(\tau A_1,\tau A_2) ( F_{ki}  -g(t_{ki})) - y(t_{n+1})+y_{n+1} \notag \\
&\quad - \tau \sum_{k=0}^{n} \mathrm{e}^{(n-k)\tau A}  \sum_{i=1}^{2} \big(\widetilde{b}_i (\tau (A_1+A_2))-b_i (\tau A_1,\tau A_2) \big) g(t_{ki}), \label{Eqn:split-en-recursion} \\
E_{n2} & =   \tau \sum_{k=0}^{n-1}\mathrm{e}^{(n+c_2-k-1)\tau A }\sum_{m=1}^2  {b}_m(\tau A_1,\tau A_2)(F_{km}-g(t_{km}))+ \tau  {a}_{21}(\tau A_1,\tau A_2))(F_{n1}-g(t_{n1}))   \notag \\
&\quad 
- y(t_{n2})+Y_{n2} -\tau \sum_{k=0}^{n-1} \mathrm{e}^{(n-k)\tau A}  \sum_{m=1}^{2} \big(\widetilde{b}_m (\tau (A_1+A_2))-b_m (\tau A_1,\tau A_2) \big) g(t_{km}) \notag \\
&\quad  -\tau  \big( \widetilde{a}_{21} (\tau (A_1+A_2)) -a_{21} (\tau A_1,\tau A_2)  \big) g(t_{n1}). \label{Eqn:split-eni-recursion} 
\end{align}
Note that $E_{n1}=e_{n}$. The analysis of $\|E_{ni}\|_\alpha$ ($i=1,2$) is quite similar to that in the proof of Theorem \ref{Thm:semilinear-first}. We just need the estimates of the terms involving split approximations. Since $u_0\in D(A)$, from \cite[Proposition 7.1.10]{Lunardi95Book} we have $u\in L^\infty((0,T);D(A))\cap C((0,T];D(A))$, which implies $\sup_{t\in[0,T]} \|g(t)\|_{D_{A_i}(\kappa,\infty)} \leq \infty$. Noting that $\widetilde{b}_m (\tau (A_1+A_2))-b_m (\tau A_1,\tau A_2) $ are actually linear combinations of $\varphi_k(\tau(A_1+A_2))-k!\varphi_k(\tau A_1) \varphi_k(\tau A_2)$, it suffices to consider the following estimate
\begin{align}
&\Bigg\|\tau \sum_{k=0}^{n-1} \mathrm{e}^{(n-k)\tau A}  \big(\varphi_k(\tau(A_1+A_2))-k!\varphi_k(\tau A_1) \varphi_k(\tau A_2)  \big) g(t_{km})\Bigg\|_{\alpha}  \notag \\
&\quad \leq C\Bigg\|\tau^3 \sum_{k=0}^{n-3} \mathrm{e}^{(n-k)\tau A} \int_0^1\int_0^1 \int_\eta^\xi    A_1\mathrm{e}^{\tau(1-\sigma)A_1} A_2\varphi_1(\tau(1-\eta)A_2)  \,\mathrm{d}\sigma \, \mathrm{d}\eta  \,\mathrm{d}\xi \cdot g(t_{km})\Bigg\|_{\alpha}  \notag \\
&\quad \quad + \sum_{k=n-2}^{n-1}\big\|\tau \mathrm{e}^{(n-k)\tau A}  \big(\varphi_k(\tau(A_1+A_2))-k!\varphi_k(\tau A_1) \varphi_k(\tau A_2)  \big) g(t_{km}) \big\|_{\alpha} \notag\\
&\quad \leq  C\tau^3 \sum_{k=0}^{n-3} \|\mathrm{e}^{\frac{1}{2} (n-k)\tau A}\|_{L(X,D_A(\alpha,1))}    \|A_l\mathrm{e}^{\frac{1}{2} (n-k)\tau A_{l}}\| \|A_i\mathrm{e}^{\frac{1}{2} (n-k)\tau A_i}\|_{L(D_{A_i}(\kappa,\infty),X)}\|g(t_{km})\|_{D_{A_i}(\kappa,\infty)} \notag\\
&\quad\quad + C\tau \sum_{k=n-2}^{n-1} \|\varphi_k (\tau A) - k!\varphi_k (\tau A_1,\tau A_2)\|_{L(D_{A_i}(\kappa,\infty),D_A(\alpha,1))} \|g(t_{k,m})\|_{D_{A_i}(\kappa,\infty)} \notag \\
&\quad\leq C\tau^{1+\kappa-\alpha},\notag
\end{align}
where $l \in \{1,2\}\backslash \{i\}$ and Lemma \ref{Est:split-varphi} was used. In addition, we have
$$
\begin{aligned}
&\big\| \tau  \big( \widetilde{a}_{21} (\tau (A_1+A_2)) -a_{21} (\tau A_1,\tau A_2)  \big) g(t_{n1}) \big\|_{\alpha} \notag \\
&\qquad \leq C\big\| \tau  \big( \widetilde{a}_{21} (\tau (A_1+A_2)) -a_{21} (\tau A_1,\tau A_2) \big) \big\|_{L(D_{A_i}(\kappa ,\infty),D_A(\alpha,1))}  \| g(t_{n1})  \|_{D_{A_i}(\kappa ,\infty)} \notag \\
&\qquad \leq C\tau^{1+\kappa-\alpha}.
\end{aligned}
$$
By using \eqref{Est:varphikvarphik-X-DAalpha1}, we conclude that 
$$
\max_{j=1,2}\|E_{nj}\|_\alpha \leq C\tau^{1+\kappa-\alpha}(\left|\log \tau \right| +1), \quad 0\leq t_n <T,
$$
and we obtain from \eqref{Eqn:split-en-recursion} that 
$$
\|u_n-u(t_n)\|_\alpha \leq C\tau^{1+\kappa-\alpha}(\left|\log \tau \right| +1), \quad 0\leq t_n \leq T.
$$
The desired estimate in $D_A(\gamma,1)$ for any $\gamma\in[\alpha,1)$ can be obtained in a similar manner as that presented in the proof of Theorem \ref{Thm:semilinear-first}. The proof is completed.
\end{proof}

We conclude this section by revisiting the examples discussed in Section \ref{Sec:example}. We focus on the nonsmooth data error estimates for the split exponential integrators.

Let $\Omega=(0,1)^2$. We first note that \eqref{Eqn:DAinfty} %and \eqref{Eqn:DAinfty-12} 
also holds for $\Omega$ being a rectangle. This follows from standard arguments, using the Seeley extension of $C^k$ functions. We split the operator $A$ as $A=A_1+A_2$, where the operators $A_1$ and $A_2$, defined by
$$
\begin{aligned}
&A_1:D(A_1)\subset C(\overline{\Omega}) \to C(\overline{\Omega}): u \mapsto \varepsilon^2 \partial_{xx} u, \\
&A_2:D(A_2)\subset C(\overline{\Omega}) \to C(\overline{\Omega}): u \mapsto \varepsilon^2 \partial_{yy} u, 
\end{aligned}
$$
are sectorial with domains
$$
\begin{aligned}
&D(A_1) =  \{ u \in C(\overline{\Omega}) :  \partial_{xx} u \in C(\overline{\Omega}),~u(0,y)=u(1,y)=0~\mbox{for}~y\in [0,1]\}, \\
&D(A_2) =  \{ u \in C(\overline{\Omega}) :  \partial_{yy} u \in C(\overline{\Omega}),~u(x,0)=u(x,1)=0~\mbox{for}~x\in [0,1]\}.
\end{aligned}
$$

First consider the numerical solution of the Allen-Cahn equation \eqref{Eqn:allen-cahn} by the split version of the exponential Euler method \eqref{Eqn:split-euler}. If $u_0\in C(\overline{\Omega})$, the error for $\tau$ sufficiently small is bounded by
$$
\|u_n-u(t_n)\|_{C(\overline{\Omega})} \leq C\tau(\left|\log \tau \right| +1).
$$
If further $u_0 \in C_0^{2\gamma}(\overline{\Omega})$ for some $\gamma\in (0,1)$, it holds 
$$
  \|g(t)\|_{D_{A_i}(\gamma,\infty)} \leq C  \|g(t)\|_{C_0^{2\gamma}(\overline{\Omega})}\leq C,\quad 0\leq t \leq T,
$$ 
where the embedding $C_0^{2\gamma}(\overline{\Omega}) \subset D_{A_{i}}(\gamma,\infty)$ was used. This embedding can be justified as follows. Define the operator $A_x$ by $A_x w = \partial_{xx} w$ for $w\in D(A_x)=\{w\in C^2([0,1]): w(0)=w(1)=0\}$.  For $v\in C_0^{2\gamma}(\overline{\Omega})$, it holds
$$
\begin{aligned}
\|v\|_{D_{A_{1}}(\gamma,\infty)}& =\|v\|_{C(\overline{\Omega})}+
\sup_{t\in (0,1)} \|t^{1-\gamma} A_1 \mathrm{e}^{t A_1} v \|_{C({\overline{\Omega}})}  \\
& =\|v\|_{C(\overline{\Omega})}+ \sup_{t\in (0,1)}  \max_{y\in [0,1]} \max_{x\in[0,1]} |t^{1-\gamma} (A_1 \mathrm{e}^{t A_1} v)(x,y) | \\
&=  \|v\|_{C(\overline{\Omega})} + \sup_{t\in (0,1)} \max_{y\in [0,1]}  \|t^{1-\gamma} A_{x}  \mathrm{e}^{t A_x} v(\cdot,y) \|_{C([0,1])}  \\
&\leq  \|v\|_{C(\overline{\Omega})} +\max_{y\in [0,1]} \|v(\cdot,y)\|_{D_{A_x}(\gamma,\infty)} \\
&\leq  \|v\|_{C(\overline{\Omega})} +C\max_{y\in [0,1]} \|v(\cdot,y)\|_{C_0^{2\gamma}([0,1])} \leq C\|v\|_{C_0^{2\gamma}(\overline{\Omega})} .
\end{aligned}
$$
Therefore, the results of Theorem \ref{Thm:semilinear-split-first} can be sharpened to yield the following estimate
$$
\|u_n-u(t_n)\|_{C_0^{2\gamma}(\overline{\Omega})} \leq C\tau(\left|\log \tau \right|+1).
$$
If further the compatibility condition $Au_0 \in C(\overline{\Omega})$ holds, then for any $\gamma\in [0,1)$ the error of the numerical solution of \eqref{Eqn:allen-cahn} by the split integrator ERK2L \eqref{ERK-split} is bounded by
$$
\|u_n-u(t_n)\|_{C_0^{2\gamma}(\overline{\Omega})} \leq C\tau^{2-\epsilon-\gamma}(\left|\log \tau \right| +1),
$$
where $\epsilon>0$ is arbitrarily small. This is a direct application of Theorem \ref{Thm:semilinear-split-first}, using the fact that $f:D(A)\to D_{A_i}(\kappa,\infty)$ is continuous for any $\kappa \in [0,1)$.

Now, we consider the numerical solution of the Burgers' equation \eqref{Eqn:Burgers} by the split integrators. The operators $A_1$ and $A_2$ can be defined in the same manner as before, except that $\varepsilon^2$ is replaced with $\nu$. The domains of the operators remain the same. 
If $u_0\in C_0^{2\gamma}(\overline{\Omega})$ for some $\gamma\in [\tfrac12,1)$, the error of the split version of the exponential Euler method \eqref{Eqn:split-euler} for $\tau$ sufficiently small is bounded by
$$
\|u_n-u(t_n)\|_{C_0^{2\gamma}(\overline{\Omega})} \leq C\tau^{1/2}(\left|\log \tau \right| +1).
$$
This estimate provides a sharper bound compared to Theorem \ref{Thm:semilinear-split-second}, by exploiting the continuity of $f:C_0^{2\gamma}(\overline{\Omega})\to C_0^{2\gamma-1}(\overline{\Omega})$.
Consider the numerical solution of \eqref{Eqn:Burgers} by the split integrator ERK2L \eqref{ERK-split}.  
If further $Au_0+f(u_0)\in C_0^1(\overline{\Omega})$, then we can verify that 
$$
\|u(t)\|_{C_0^{3-2\epsilon}(\overline{\Omega})} \leq C,\quad 0\leq t\leq T,\quad
\mbox{and}\quad 
\|g(t)\|_{D_{A_i}(1-\epsilon,\infty)} \leq C,\quad 0\leq t\leq T,
$$
where $\epsilon>0$ is arbitrarily small. Therefore, for any $\gamma \in[\tfrac12,1)$ the error is bounded by
$$
\|u_n-u(t_n)\|_{C_0^{2\gamma}(\overline{\Omega})} \leq C\tau^{2-\epsilon-\gamma}(\left|\log \tau \right| +1).
$$

\section*{Acknowlegments}
Qiumei Huang is supported by the National Natural Science Foundation of China (No.~12371385). 
Gangfan Zhong is supported by the China Scholarship Council (CSC) joint Ph.D. student
scholarship (Grant 202406540082).

%\section*{Acknowlegments}
%Qiumei Huang is supported by the National Natural Science Foundation of China (No.~12371385). 
%Gangfan Zhong is supported by the China Scholarship Council (CSC) joint Ph.D. student
%scholarship (Grant 202406540082).

%%%%%%%%%%%%%%%%%%%%
\bibliographystyle{abbrv}

%% 如下给出 Bib 数据库的路径和名称
\bibliography{/Volumes/Bibliography/gfzhong_bib/Mathematics.bib}

%\bibliography{/Users/zhongliuqiang/Math/Bibtex/0Student/Math}
%\bibliography{/Volumes/Book/数据库/Math_2021_1222.bib}

\end{document}